\documentclass[final]{siamltex}
\usepackage{amssymb}
\usepackage{mathrsfs}
\usepackage{amsfonts, dsfont}
\usepackage{graphicx}
\usepackage{colortbl,dcolumn}
\usepackage{amsmath}
\usepackage{psfrag}
\usepackage{booktabs}
\usepackage{cite}
\usepackage{url}

\numberwithin{equation}{section}
\usepackage[top=1.4in, bottom=1.4in, left=1.2in, right=1.2in]{geometry}

\newcommand{\R}{\mathbb{R}}
\newcommand{\N}{\mathbb{N}}

\newcommand{\dd}{\text{d}}

\newtheorem{thm}{Theorem}[section]

\newtheorem{assumption}[thm]{Assumption}
\allowdisplaybreaks[1]

\title{Error estimates of semi-discrete and fully discrete finite element methods for the Cahn-Hilliard-Cook equation
\thanks
{
R.Q. was supported by National Science Foundation of China (11701073).  
X.W. was supported by National Science Foundation of China (11671405, 11971488, 91630312)
and Natural Science Foundation of Hunan Province for Distinguished Young Scholars (2020JJ2040).
Dr. Meng Cai is gratefully acknowledged for bringing a few typos into our notice.
Also, the authors want to thank the Tianyuan Mathematical Center in Northeast China
for the hospitality and Prof. Xiaobing Feng from University of Tennessee for his helpful comments
when this work was first presented in a conference in June of 2018, hosted by the center.
}
}

\author{
Ruisheng Qi\thanks{School of Mathematics and Statistics,
       Northeastern University at Qinhuangdao, Qinhuangdao, China
({\tt qirsh@neuq.edu.cn}).}
\and Xiaojie Wang\thanks{School of Mathematics and Statistics,
       Central South University, Changsha 410083, Hunan, China
({\tt x.j.wang7@csu.edu.cn, x.j.wang7@gmail.com}).}
}

\begin{document}

\maketitle

\begin{abstract}
In two recent  publications [Kov{\'a}cs, Larsson, and Mesforush, SIAM J. Numer. Anal. 49(6), 2407-2429, 2011]
and [Furihata, et al., SIAM J. Numer. Anal. 56(2), 708-731, 2018],
strong convergence of the semi-discrete and fully discrete finite element methods is,
respectively, proved for the Cahn-Hilliard-Cook (CHC) equation, but without convergence rates revealed.
The present work aims to fill the left gap, by recovering strong convergence rates of (fully discrete) finite element methods
for the CHC equation. More accurately, strong convergence rates of a full discretization are obtained,  based on Galerkin
finite element methods for the spatial discretization and the backward Euler method for the temporal discretization. 
It turns out that the convergence rates heavily depend on 
the spatial regularity of the noise process. Different from the stochastic Allen-Cahn equation, 
the presence of the unbounded elliptic operator in front of the cubic nonlinearity
in the underlying model makes the error analysis much more challenging and demanding.
To address such difficulties, several new techniques and error estimates are developed.
Numerical examples are finally provided to confirm the previous findings.
%
%
\end{abstract}




\section{Introduction}
 Over the last  twenty years,   numerical approximations of stochastic partial
differential equations (SPDEs) with globally Lipschitz coefficients have been extensively and well studied,
see the monographs \cite{JK11,kruse2014strong,lord2014introduction} and references therein. By contrast,
numerical analysis of SPDEs with non-globally Lipschitz coefficients is, in our opinion, at an early stage and far from being well-understood.
A typical SPDE model with non-globally Lipschitz coefficients is the stochastic Allen-Cahn equation, which has been
numerically studied by many researchers recently, see, e.g.,  \cite{becker2017strong,becker2019strong,brehier2018strong,brehier2018analysis,brehier2018weak,feng2017finite,liu2018strong,liu2018strong-multiplicative,cui2018weak,cui2019strong,feng2018strong,gyongy2016convergence,jentzen2015strong,jentzen2018exponential,kovacs2015discretisation,kovacs15backward,Majee2017optimal,campbell2018adaptive}.
As another prominent SPDE model with non-globally Lipschitz coefficients,
the Cahn-Hilliard-Cook (CHC) equations, also known as the stochastic Cahn-Hilliard equation in some literature, 
are, however, much less investigated.
As far as we know, only a few publications are devoted to numerical studies of the CHC equation
\cite{cardon2000implicit,cui2018strong,hutzenthaler2014perturbation,kovacs2011finite,furihata2018strong,
li2016unconditionally}.
Particularly,  strong convergence of the semi-discrete and fully discrete finite element methods is,
respectively, proved in \cite{kovacs2011finite} and \cite{furihata2018strong} for the CHC equation,
but without convergence rates recovered.
The present article attempts to fill the left gap, by recovering strong convergence rates of the (fully discrete) finite element methods
for the CHC equation.

Let $D\subset \mathbb{R}^d, d\in\{1,2,3\}$ be a bounded convex spatial domain with smooth boundary and
 let
 $H:=L_2(D; \mathbb{R})$
  be the real separable Hilbert space
 endowed with the usual inner product  and  norm.
 Throughout the paper we are interested in the following Cahn-Hilliard-Cook equation perturbed by noise  in $\dot{H}:=\{v\in H:\int_D v\,dx=0\}$,
 \begin{align}\label{eq:stochastic-allen-Hilliard-equation}
\begin{split}
\left\{
\begin{array}{ll}
\dd u - \Delta w \, \dd t= \dd W,&\; \text{ in }\; D\times (0,T],
\\
w=- \Delta u+f(u), &\; \text{ in } \; D\times (0,T],
\\
\frac{\partial u}{\partial n}=\frac{\partial w}{\partial n}=0,& \; \text{ in } \;\partial D\times(0,T],
\\
u(0,x)=u_0,&\; \text{ in } \; D,
\end{array}
\right.
\end{split}
\end{align}
where $\Delta = \sum_{ i = 1 }^{d} \tfrac{ \partial ^2}{ \partial x_i^2} $,  $f(s)=s^3-s, s\in \mathbb{R}$
and $ \frac{\partial }{\partial n} $ denotes the outward normal derivative on $ \partial D $. 
Following the framework of \cite{da1996stochastic} we rewrite \eqref{eq:stochastic-allen-Hilliard-equation} as
an abstract evolution equation of the form,
\begin{align}\label{eq:abstract-SCHE}
\left\{\begin{array}{ll}
\dd X(t)+ A\big(AX(t)+F(X(t))\big)\,\dd t=\,\dd W(t),\;&t\in(0,T],
\\
X(0)=X_0,
\end{array}
\right.
\end{align}
where $-A$ is the Laplacian with homogeneous Neumann boundary conditions and $- A^2$
 generates an analytic semigroup $E(t)$ in $\dot{H}$. Similarly as in \cite{kovacs2011finite,furihata2018strong},
 $\{W(t)\}_{t\geq 0}$ is assumed to be an $\dot{H}$-valued $Q$-Wiener process on a filtered probability  space
 $(\Omega, \mathcal{F},  \mathbb{P}, \{\mathcal{F}_t\}_{t\geq 0})$.
  The nonlinear mapping $F$ is supposed to be a Nemytskij operator, given by $F(u)(x)=f(u(x))$, $x\in D$.

The deterministic version of such equation is used to
 describe the complicated phase separation and coarsening phenomena in a melted alloy \cite{cahn1961spinodal,cahn1958free,cahn1971spinodal} that is
quenched to a temperature at which only two different concentration phases can exist stably.
For such model, $u$ represents the concentration of an alloy and $w$ models the chemical potential.
The corresponding numerical study, e.g., can be consulted in \cite{elliott1992error}.
 Concerning the stochastic version,  Da Prato and Debussche \cite{da1996stochastic} have already proved
 the existence and uniqueness of the solution to \eqref{eq:abstract-SCHE}.
 The space-time regularities of the weak solution of   \eqref{eq:abstract-SCHE}
 have been further examined in \cite{larsson2011finite,furihata2018strong}.   As the first goal, this work aim to provide improved 
 regularity results for the solution  to  \eqref{eq:abstract-SCHE} based on existing ones from \cite{da1996stochastic,larsson2011finite,furihata2018strong}.
    Under further assumptions
   specified later, particularly including
   \begin{align}\label{eq:assumption-on-Q}
   \|A^{\frac{\gamma-2}2}Q^{\frac12}\|_{\mathcal{L}_2}<\infty,\; \text{ for some }\;\gamma\in[3,4],
   \end{align}
Theorems \ref{theorem:preliminary-theorem}, \ref{them:regulairty-mild-solution} assert that, 
\eqref{eq:abstract-SCHE} admits a  unique mild solution $X(t), t \in [0, T]$, given by
  \begin{align}\label{them:eq-mild-solution-stochastic-equation}
X(t)
=
E(t)X_0
-
\int_0^tE ( t - s ) A F(X(s))\,\mathrm{d} s
+
\int_0^tE(t-s)\mathrm{d} W(s),
\end{align}
which   enjoys the following spatial-temporal regularity properties,
  \begin{align*}
  X\in L_\infty([0,T];L^p(\Omega; \dot{H}^\gamma)),\;\forall p\geq 1,
  \end{align*}
  and
  for $\forall p\geq 1$ and $0\leq s<t\leq T$,
\begin{align*}
\|X(t)-X(s)\|_{L^p(\Omega, \dot{H}^\beta)}
\leq
C(t-s)^{\min\{\frac12,\frac{\gamma-\beta}4\}},\;\beta\in[0,\gamma].
\end{align*}
Here
$\dot{H}^\alpha:=\text{dom}(A^{\frac\alpha2}), \alpha\in \mathbb{R}$
and the parameter $\gamma\in[3,4]$ coming from \eqref{eq:assumption-on-Q} quantifies the spatial regularity of
the covariance operator $Q$ of the driving noise process.

The second aim of this article is to derive error estimates for finite element approximations of the stochastic problem \eqref{eq:abstract-SCHE}.
By $\dot{V}_h \subset H^1(D)\cap \dot{H}$  we denote the space of continuous functions that are piecewise polynomials
of degree at most $r-1$,  for $r\in\{2,3,4\}$ in dimension $d=1$ and $r=2$ in dimension $d\in\{2,3\}$, 
and by $X_h(t)\in \dot{V}_h$ the finite element spatial approximation of the mild solution $X$, represented by
\begin{align} \label{eq:intro-Xh}
X_h(t)=E_h(t)P_hX_0
-
\int_0^tE_h(t-s) A_h P_h F(X_h(s))\,\dd s
+
\int_0^tE_h(t-s)P_h\,\dd W(s),
\quad
t\in[0,T].
\end{align}
Here $ h > 0$ is the mesh size and $E_h(t):=e^{-tA^2_h}$ is the strongly continuous semigroup generated by
the discrete operator $ - A_h^2$. The resulting spatial approximation error, as implied by
Theorem \ref{them:error-estimates-semi-problem}, is measured as follows,
\begin{align}\label{eq:error-semidiscrete-problem}
\|X(t)-X_h(t)\|_{L^p(\Omega;\dot{H})}
\leq
Ch^\kappa|\ln h|,
\quad
\kappa := \min\{\gamma,r\}.
\end{align}
To arrive at it, we introduce an auxiliary approximation process $ \widetilde{X}_h $, defined by
\begin{align}\label{eq:mild-solution-semi-auxiliary-problem}
\widetilde{X}_h(t)=E_h(t)P_hX_0
-
\int_0^tE_h(t-s) A_h P_h F(X(s))\,\dd s
+
\int_0^tE_h(t-s)P_h\,\dd W(s),
\quad
t\in[0,T],
\end{align}
and split the considered error  $ \|X(t)-X_h(t)\|_{L^p(\Omega;\dot{H})} $ into two parts:
\begin{equation}
\|X(t)-X_h(t)\|_{L^p(\Omega;\dot{H})}
\leq
\|X(t) - \widetilde{X}_h(t) \|_{L^p(\Omega;\dot{H})}
+
\|\widetilde{X}_h(t) - X_h(t)\|_{L^p(\Omega;\dot{H})}.
\end{equation}
In a semigroup framework, one can straightforwardly treat the first error term
and show $\|X(t)-\widetilde{X}_h(t)\|_{L^p(\Omega;\dot{H})}=O(h^\kappa|\ln h|)$, with the aid of the well-known
estimates for the error operators $\Psi_h(t):=E(t)-E_h(t)P_h$ and $\Phi_h(t):=E(t)A-E_h(t)A_hP_h$ and uniform
 moment bounds of  $\widetilde{X}_h(t)$ and  $X(t)$.
 Further, we subtract \eqref{eq:intro-Xh} from
\eqref{eq:mild-solution-semi-auxiliary-problem} to eliminate the stochastic convolution
and the remaining term $\widetilde{e}(t):=\widetilde{X}_h(t)-X_h(t)$ satisfies
\begin{align}
\label{eq:intro-error-hat-PDE}
\dd\widetilde{e}_h(t)
+
A_h^2\widetilde{e}_h(t) \dd t
=
A_hP_h \big( F(X_h(t)) - F(X(t)) \big)\dd t,
\quad
\widetilde{e}_h(0)=0,
\end{align}
whose solution is given by
\begin{align}\label{eq:mild-solution-e_h}
\widetilde{e}_h(t)
=
\int_0^tE_h(t-s)A_hP_h \big(F(X_h(s))-F(X(s))\big)\,\dd s.
\end{align}
Note that the tough term $\| \widetilde{e}( t ) \|_{L^p(\Omega;\dot{H})}$ can not be handled directly
due to the presence of $A_h$ before the nonlinearity.
 However, we turn things around and derive
 $\big\|\int_0^t|\widetilde{e}_h(s)|_1^2\,\dd s\big\|_{L^p(\Omega; \mathbb{R})} = O (h^{2\kappa} | \ln h|^2 )$ instead,
 after fully exploiting \eqref{eq:intro-error-hat-PDE},  the monotonicity of the nonlinearity, regularity
 properties of $X_h(t)$, $\widetilde{X}_h(t)$ and $X(t)$,  and the previous error estimate for $\|X(t)-\widetilde{X}_h(t)\|_{L^p(\Omega;\dot{H})}$.
  Equipped with the key error estimate of $\big\|\int_0^t|\widetilde{e}_h(s)|_1^2\,\dd s\big\|_{L^p(\Omega; \mathbb{R})}$
  and \eqref{eq:mild-solution-e_h},  we can smoothly
  show $\|\widetilde{X}_h(t)-X_h(t)\|_{L^p(\Omega;\dot{H})}=O(h^\kappa| \ln h|)$
 (see \eqref{eq:error-widetild-e-semi}-\eqref{eq:bound-L2-term}) and therefore obtain \eqref{eq:error-semidiscrete-problem}.

Let $k=T/N$, $N\in \mathbb{N}$ be a uniform time step-size.
After discretizing the stochastic problem \eqref{eq:abstract-SCHE} by  the finite element method in space and the backward Euler  scheme in time, we also investigate the resulting fully discrete scheme, given by
\begin{align*}
X^n_h
=
E_{k,h}  X^{n-1}_h
-
k E_{k,h} A_h  P_h  F ( X_h^{ n } )
+
E_{k,h} P_h \Delta W_{n},
\end{align*}
where $E_{k,h}:=(I+kA_h^2)^{-1}$ and $X_h^n$ is regarded as the fully discrete approximation of $X(t_n)$.
By essentially exploiting discrete analogue of arguments as used in the semi-discrete case,
one can obtain the following strong approximation error bound
\begin{align*}
\|X(t_n)-X_h^n\|_{L^p(\Omega;\dot{H})}
\leq
C(h^\kappa| \ln h|+k^{\frac\kappa4}| \ln k|),
\quad
\kappa := \min\{\gamma,r\}.
\end{align*}

It is important to mention that,
the presence of the unbounded operator $A$ in front of the non-globally Lipschitz (cubic) nonlinearity in
the underlying model causes essential difficulties in the error analysis for the approximations
 and the error analysis becomes much more challenging than that of the stochastic Allen-Cahn equation
 (see \cite{qi2018optimal} and relevant comments in \cite{kovacs2011finite,furihata2018strong}).
More specifically, our error analysis heavily relies on the new approach mentioned before,
a priori strong moment bounds of the numerical approximations,
and a variety of error estimates for the finite element approximation of  the corresponding deterministic linear problem.
Some estimates can be derived from existing ones in \cite{kovacs2014erratum, larsson2011finite,furihata2018strong}.
Nevertheless, estimates available in \cite{kovacs2014erratum, larsson2011finite,furihata2018strong} are
far from being enough for the purpose of the error analysis. For example, the
strong moment bounds \eqref{lem:eq-bound-solution-semi} and \eqref{lem:eq-bound-solution-full-stochatic} and
the error estimates of integral form such as \eqref{lem:eq-chemical-potial-integrand}, \eqref{lem:eq-chemical-potial-integrand-II},
\eqref{lem:error-deterministc-potial-full-integrand} and \eqref{lem:error-deterministc-potial-full-integrand-II} are completely new.

Finally, we add some comments on a few closely relevant works.  A finite  difference scheme was examined in \cite{cardon2000implicit} 
for the problem \eqref{eq:abstract-SCHE} and convergence in probability was
established with rates. Hutzenthaler and  Jentzen   \cite{hutzenthaler2014perturbation} used a general perturbation theory
and exponential integrability properties of the exact and numerical solutions
 to prove strong convergence rates for the spatial spectral Galerkin approximation (no time discretizaton)
  in one spatial dimension.  In \cite{larsson2011finite,furihata2018strong},
  strong convergence of finite element methods for \eqref{eq:abstract-SCHE} was proved, but with no rate obtained.
   The analysis in \cite{larsson2011finite,furihata2018strong} is based on proving a priori moment bounds with large exponents and
   in higher order norms using energy arguments and bootstrapping followed by a pathwise Gronwall argument in the mild solution setting.
 %
 %
    Before submitting the early version of the present work to arXiv in late December of 2018,
    we were also aware of an interesting preprint \cite{cui2018strong} submitted to arXiv in mid-December of 2018.
    There strong convergence rates of a fully discrete scheme are obtained,
    done by a spatial spectral Galerkin method and a temporal accelerated implicit Euler method for the CHC equations.
    To the best of our knowledge,   strong convergence rates of finite element methods for the CHC equations
    are missing in the existing literature and this article fills the left gap.

  The outline of this paper is as follows. In the next section, some preliminaries are collected and certain assumptions 
  are made to ensure well-posedness of the considered problem. 
  Section 3 is devoted to the uniform moment bounds of the semi-discrete
  finite element approximation. Based on the uniform moment bounds obtained in section 3, we derive the error estimates
  for the semi-discrete problem in section 4. Section 5 focuses on the uniform moment bounds of the fully discrete
  approximations and section 6 provides error estimates of the backward Euler-finite element full discretization.
  In section \ref{sect:numerical-section}, numerical examples are provided to confirm the previous findings.

\section{The Cahn-Hilliard-Cook equation}
Throughout this paper,  we use $\N$ to denote the set of all positive integers and denote $\N_0 = \{ 0 \} \cup \N$.
Given a separable $\mathbb{R}$-Hilbert space $\big(H, \;\big(\cdot,\cdot\big), \;\|\cdot\|\big)$,
 by $\mathcal{L}(H)$ we denote the Banach space of all linear bounded operators from $H$ to $H$. Also, we denote  by $\mathcal{L}_2(H)$
the Hilbert space consisting of  all Hilbert-Schmidt operators from
$H$ into $H$, equipped with the inner product and the norm,
\begin{align*}
\big<\Gamma_1,\Gamma_2\big>_{\mathcal{L}_2(H)}
=
\sum_{j=1}^\infty \big<\Gamma_1 \phi_j,\Gamma_2\phi_j\big>,
\qquad
\|\Gamma\|_{\mathcal{L}_2(H)}
=
\bigg(
\sum_{j=1}^\infty \big\|\Gamma \phi_j\|^2
\bigg)^\frac12
,
\end{align*}
independent of the choice of orthonormal basis $\{\phi_j\}$ of $H$. If $\Gamma\in\mathcal{ L}_2(H)$ and $L\in \mathcal{L}(H)$, then $\Gamma L, L\Gamma\in \mathcal{L}_2(H)$ and
\begin{align*}
\|\Gamma L\|_{\mathcal{L}_2(H)}
\leq
\|L\|_{\mathcal{L}(H)}
\|\Gamma\|_{\mathcal{L}_2(H)}
,
\quad
\|L\Gamma \|_{\mathcal{L}_2(H)}
\leq
\|L\|_{\mathcal{L}(H)}
\|\Gamma\|_{\mathcal{L}_2(H)}.
\end{align*}


\subsection {Abstract framework and main assumptions}
In this subsection, we formulate main assumptions concerning the operator $A$,
the nonlinear mapping, the noise process  and the initial data.

\begin{assumption}\label{assum:linear-operator-A} (Linear operator $A$)
Let $D$ be a bounded convex domain in $\mathbb{R}^d$ for $d\in\{1,2,3\}$ with a sufficiently smooth boundary $\partial D$ and
let $H=L_2(D; \mathbb{R})$  be the real separable Hilbert space endowed
 with the usual inner product $(\cdot,\cdot)$ and the associated norm $\|\cdot\|=(\cdot,\cdot)^{\frac12}$.
 Let $\dot{H} = \big\{v \in H : (v , 1) = \int_Dv\,\mathrm{d}x = 0\big\}$ and let $ - A = \Delta $ be the Neumann Laplacian,
 with the domain $ dom(A) := \big\{v\in H^2(D)  \cap \dot{H}  \colon \frac{\partial v}{\partial n}=0, \;on\;\partial D\big\}$.
\end{assumption}

Such assumptions guarantee that the operator $A$ is positive definite, self-adjoint, unbounded, linear on $\dot{H}$ with
compact inverse.  Let
$
P \colon H  \rightarrow \dot{H}
$
denote a generalized orthogonal projection, given by $ P v = v - |D|^{-1}\int_D v\,\dd x$.
Then $(I-P)v = |D|^{-1}\int_D v\,\dd x$ is the average of $v$ and it is not difficult to check
\begin{equation}
\label{eq:P-Lq-bound}
\| Pv \|_{L_q} \leq 2 \| v \|_{L_q},
\quad
q \geq 2.
\end{equation}
Here and below, by $L_r ( D;\mathbb{R} ), r\geq 1$ ($L_r(D)$ or $L_r$ for short) we denote a Banach space
consisting of $r$-times integrable functions.
When extended to $H$ as $A v:= APv$, for $v\in H$, the linear operator $A$ has
an orthonormal basis  $\{e_j\}_{j \in \N_0 }$ of $H$ with corresponding eigenvalues   $\{\lambda_j\}_{j \in \N_0}$   such that
\begin{align}\label{eq:eigenvalue-A-H}
0=\lambda_0<\lambda_1\leq\lambda_2\leq\cdots\leq\lambda_j\leq\cdots, \quad\lambda_j\rightarrow\infty .
\end{align}
Note that the first eigenfunction is a constant, i.e., $e_0=|D|^{-\frac12}$ and $\{e_j\}_{j \in \N}$ forms an orthonormal basis of $\dot{H}$.
By the spectral theory, we can define the fractional powers
of $A$  on $\dot{H}$ in a simple way, e.g., $A^\alpha v=\sum_{j=1}^\infty\lambda_j^\alpha\left(v,e_j\right)e_j$, $\alpha\in \mathbb{R}$.
Define
the inner product $\left( \cdot,  \cdot\right)_{\alpha}$ and the associated norm
 $|\cdot|_\alpha:=\|A^{\frac\alpha 2}\cdot\|$, given by
 \begin{align*}
|v |_{\alpha}=\|A ^{\frac\alpha2} v \|
=
\Big(\sum_{j=1}^{ \infty }
\lambda_{ j }^\alpha| (v, e_ { j })|^2\Big)^{\frac12},
\quad
\left( v,  w \right)_{\alpha}
=
\sum_{ j = 1}^{ \infty }
\lambda_j^\alpha
 (v, e_ { j })  (w, e_ { j }),
\quad
\alpha\in \mathbb{R}.
\end{align*}
%
Then we define the following function spaces
\begin{equation*}
\dot{H}^\alpha:= \text{dom} (A^{\frac\alpha2})
= \{ v \in \dot{ H } \colon | v |_{\alpha} < \infty \}
,
\quad
\alpha \geq 0.
\end{equation*}
Then $\dot{H}^0 = \dot{ H}$.
For negative order $ - \alpha < 0$ we define $ \dot{ H }^{-\alpha} $ by taking the closure of $\dot{H}$ with respect to $| \cdot |_{ - \alpha}$.
 It is known that for  integer $\alpha\geq 0$, $\dot{H}^\alpha$ is a subspace
of $H^\alpha(D) \cap \dot{H}$ characterized by certain boundary conditions. Additionally, 
the norm $|\cdot|_\alpha$ is  equivalent on $\dot{H}^\alpha$ to the standard Sobolev norm 
$\|\cdot\|_{H^\alpha(D)}$ for $\alpha=1,2$.
%

Thanks to  \eqref{eq:eigenvalue-A-H}, the operator $-A^2$ can
generate an analytic semigroup $E(t)=e^{-tA^2}$ on $H$ and
\begin{align}\label{eq:definition-semigroup-E(t)}
\begin{split}
E(t)v=e^{-tA^2}v
=
\sum_{j=1}^\infty e^{-t \lambda_j^2}\big(v,e_j\big)e_j
+
\big(v,e_0\big) e_0
=Pe^{-t A^2}v
+
(I-P)v,
\quad
v \in H.
\end{split}
\end{align}
By expansion in terms of the eigenbasis of $A$ and using the Parseval identity,  one can easily obtain
\begin{align}
\|A^\mu E(t)\|_{\mathcal{L}(\dot{H})}
&
\leq
 Ct^{-\frac\mu2},\; t>0,\; \mu\geq 0,
 \label{I-spatial-temporal-S(t)}
 \\
 \|A^{-\nu}(I-E(t))\|_{\mathcal{L}(\dot{H})}
 &
 \leq
 Ct^{\frac\nu2},\quad t\geq0,\;\nu\in[0,2],
 \label{II-spatial-temporal-S(t)}
 \\
\int_{\tau_1}^{\tau_2} \|A^\varrho E(s) v\|^2\,\dd s
&
\leq
C ( \tau_2 - \tau_1 ) ^{1-\varrho}\|v\|^2,\;\forall v\in \dot{H}, \varrho\in[0,1],
\tau_2 \geq \tau_1 \geq 0, 
\label{III-spatial-temporal-S(t)}
\\
 \Big\|A^{2\rho}\int_{\tau_1}^{\tau_2}
         E(\tau_2-\sigma)v\,\dd \sigma\Big\|
 &
 \leq
 C ( \tau_2-\tau_1 ) ^{1-\rho}\|v\|,\;\forall v\in \dot{H},\; \rho\in[0,1],
 \tau_2 \geq \tau_1 \geq 0.
 \label{IV-spatial-temporal-S(t)}
 \end{align}
Throughout the paper, $C$ denotes a generic positive constant that may change from line to line.
The next assumption specifies the nonlinearity of the considered equation.
\begin{assumption}\label{assum:nonlinearity}(Nonlinearity)
Let $F \colon L_6(D ; \mathbb{R}) \rightarrow H$ be a deterministic mapping given by
\begin{align*}
F(v)(x)=f(v(x))=v^3(x)-v(x),
\quad
x\in D, v\in L_6(D;\mathbb{R}).
\end{align*}
\end{assumption}
It is easy to check that, for any $v, \psi, \psi_1, \psi_2\in L_6(D)$,
\begin{align}\label{eq:definition-F-derivation}
\begin{split}
(F'(v)(\psi)\big)(x)&=f'(v(x))\psi(x)=\big(3v^2(x)-1\big)\psi(x),
\quad
x\in D,
\\
\big(F''(v)(\psi_1,\psi_2)\big)(x)&=f''(v(x))\psi_1(x)\psi_2(x)=6v(x)\psi_1(x)\psi_2(x),
\quad
x\in D.
\end{split}
\end{align}
Moreover, there exists a constant $C$ such that
\begin{align}
-(F(u)-F(v),u-v)&\leq
\|u-v\|^2,\qquad\quad u, v\in L_6(D),
\label{eq:one-side-condition}
\\
\|F(u)-F(v)\|
\leq
C \| u - v \|
( 1+&\|u\|_V^2 + \|v\|_V^2 ),
\qquad
u, v\in V,
\label{eq:local-condition}
\end{align}
where by $ V:=C(D ; \mathbb{R}) $ we denote a Banach space of continuous functions
with a usual norm.
In order to make the solution $X$ preserve the total mass, that is, $ ( I - P ) X ( t ) = ( I - P ) X_0 $,
we assume the average of the Wiener process to be zero.
\begin{assumption}\label{assum:noise-term}(Noise process)
Let $\{W(t)\}_{t\in[0,T]}$ be a standard $\dot{H}$-valued $Q$-Wiener process on the stochastic basis
$(\Omega, \mathcal{F}, \mathbb{P}, \{\mathcal{F}_t\}_{t\in[0,T]})$, where the covariance operator
 $Q\in \mathcal{L}(\dot{H})$ is bounded, self-adjoint and positive semi-definite, satisfying
\begin{align}\label{eq:ass-AQ-condition}
\|A^{\frac{\gamma-2}2}Q^{\frac12}\|_{\mathcal{L}_2}<\infty,\;for\;some\;\gamma\in[3,4].
\end{align}

\end{assumption}
\begin{assumption}\label{assum:intial-value-data}(Initial data)
Let $X_0:\Omega \rightarrow \dot{H}$ be $\mathcal{F}_0/\mathcal{B}(\dot{H})$-measurable and satisfy, for a sufficiently large number $p_0\in \mathbb{N}$,
\begin{align*}
\mathbf{E}[|X_0|_\gamma^{p_0}]<\infty,
\end{align*}
where $\gamma\in[3,4] $ is the parameter coming from \eqref{eq:ass-AQ-condition}.
\end{assumption}
%
%

\subsection{ Regularity results of the model}
This part is devoted to the well-posedness of the underlying problem \eqref{eq:abstract-SCHE} and
the space-time regularity properties of the mild solution. Existence, uniqueness, and regularity of
 weak  and mild solutions to \eqref{eq:abstract-SCHE} have been studied in
  \cite{da1996stochastic,kovacs2011finite}. The relevant result is stated as follows.
\begin{theorem}\label{theorem:preliminary-theorem}
If Assumptions \ref{assum:linear-operator-A}-\ref{assum:intial-value-data} are valid,  then the problem
\eqref{eq:abstract-SCHE} admits a weak solution $X(t)$, which is almost surely continuous and satisfies the equation
\begin{align}
\begin{split}
\langle X(t), v \rangle
-
\langle X_0, v \rangle
+
\int_0^t
\langle
X(s), A^2v
\rangle
&+
\langle
 F(X(s)), Av
 \rangle
\,\mathrm{d} s
\\
=
&
\langle
W(t), v
\rangle,
\:
a.s., 
\quad
\forall v\in \dot{H}^4 = \text{dom}(A^2), \, t \in [ 0, T].
\end{split}
\end{align}
In addition, the weak solution $X(t)$ is also a mild solution, given by \eqref{them:eq-mild-solution-stochastic-equation}, satisfying
\begin{align}\label{them:eq-bound-mild-solution-H1}
\sup_{ t \in[0,T]}\|X( t )\|_{L^p(\Omega;\dot{H}^1)}
<
\infty,\quad\forall p\geq 1.
\end{align}
\end{theorem}
To validate \eqref{them:eq-bound-mild-solution-H1}, one can simply adapt the proof of
\cite[Theorem 3.1]{kovacs2011finite}, where it was shown that
$\mathbf{E}[J(X(t))]+\mathbf{E}[\int_0^tJ'(X(s))\,\dd s]\leq C(t)$
by introducing the following  Lyapunov functional
\begin{align}\label{eq:definition-lyapunov-fun}
J(u) = \tfrac12\|\nabla u\|^2
+
\int_D\Phi(u)\,\dd x,\; u\in \dot{H}^1.
\end{align}
Here $\Phi(s) : = \frac14 (s^2-1)^2$ is a primitive of $f(s)=s^3-s$.
Evidently, the above estimate \eqref{them:eq-bound-mild-solution-H1} together with the embedding inequality 
$\dot{H}^1\subset L_6(D)$ suffices to ensure
\begin{align}\label{eq:bound-F(X)}
\begin{split}
\sup_{s\in[0,T]}\|F(X(s))\|_{L^p(\Omega;H)}
\leq
C\Big(1+\big(\sup_{s\in[0,T]}\|X(s)\|_{L^{3p}(\Omega;\dot{H}^1)}\big)^3\Big)
<\infty,
\end{split}
\end{align}
and
similarly
\begin{align}\label{eq:one-seoncd-derivation-f}
\sup_{s\in[0,T]}\|f'(X(s))\|_{L^p(\Omega;L_3)}
+
\sup_{s\in[0,T]}\|f''(X(s))\|_{L^p(\Omega;L_6)}
<\infty.
\end{align}
Furthermore, one can show further properties of the mild solution as follows.
\begin{theorem}\label{them:regulairty-mild-solution}
Let Assumptions \ref{assum:linear-operator-A}-\ref{assum:intial-value-data} be fulfilled. Then the unique mild
 solution \eqref{them:eq-mild-solution-stochastic-equation} enjoys the following regularity properties, 
\begin{align}\label{them:spatial-regularity-mild-solution}
\sup_{ t \in[0,T]} \| X(t) \|_{L^p(\Omega;\dot{H}^\gamma)}
&<\infty,\;\forall p\geq 1,
\\
\label{them:temporal-regularity-mild-stoch}
\|X(t)-X(s)\|_{L^p(\Omega;\dot{H}^{\beta})}
&\leq
C|t-s|^{\min\{\frac12,\frac{\gamma-\beta}4\}},\quad
\forall p\geq 1, \;
0 \leq s< t \leq T,
 \,
 \beta\in[0,\gamma],
\end{align}
where $\gamma\in[3,4]$ comes from Assumption \ref{assum:noise-term}.
\end{theorem}

Here and below $C$ is a generic positive constant that is also dependent of $\gamma, p$, $T$ and the initial data, 
but independent of the discretization parameters $h$ and $k$.
To prove the theorem, we introduce some basic inequalities.
Recall first the following embedding inequalities,
\begin{align}\label{eq:embedding-equatlity-I}
\dot{H}^1\subset L_6(D)\quad \text{and} \quad \dot{H}^\delta \subset C(D;\mathbb{R}),
\quad
\text{ for } \delta>\tfrac d2, \;
d\in\{1,2,3\}.
\end{align}
With \eqref{eq:P-Lq-bound} and \eqref{eq:embedding-equatlity-I} at hand, 
one can further derive, for any $\delta>\frac d2$ and any $x\in L_2(D)$,
\begin{align}\label{eq:embedding-equatlity-III}
\begin{split}
\|A^{-\frac\delta2} Px\|
=
\sup_{v\in \dot{H}}\frac{|(Px,A^{-\frac \delta2}v)|}{\|v\|}
\leq
\sup_{v\in \dot{H}}\frac{\|Px\|_{L_1} \|A^{-\frac \delta2}v\|_V} { \|v\| }
\leq
C\sup_{v\in \dot{H}}\frac{\|Px\|_{L_1} \|v\|}{\|v\|}
\leq
C\|x\|_{L_1}.
\end{split}
\end{align}
Similarly,  one can see that, for any $x\in L_{\frac65}(D)$,
\begin{align}\label{eq:embedding-equatlity-IIII}
\begin{split}
\|A^{-\frac12}P x\|
=
\sup_{v\in \dot{H}}\frac{|(Px,A^{-\frac 12}v)|}{\|v\|}
\leq
\sup_{v\in \dot{H}}\frac{\|Px\|_{L_{\frac65}} \|A^{-\frac 12}v\|_{L_6}}{\|v\|}
\leq
C\sup_{v\in \dot{H}}\frac{\|x\|_{L_{\frac65}} \|v\|}{\|v\|}
\leq
C\|x\|_{L_{\frac65}}.
\end{split}
\end{align}
Since  the norm $|\cdot|_2$ is  equivalent on $\dot{H}^2$ to
the standard Sobolev norm $\|\cdot\|_{H^2(D)}$ and $H^2(D)$
is an algebra,  one can find a constant $C >0$ such that, for any $f,g\in \dot{H}^2$,
\begin{align}\label{eq:algebra-properties-Hs}
\|fg\|_{H^2(D)}
\leq
C\|f\|_{H^2(D)}\|g\|_{H^2(D)}
\leq
C | f |_2 \,| g |_2 .
\end{align}
%
{\it Proof of Theorem \ref{them:regulairty-mild-solution}.}
We start by proving a preliminary spatial-temporal regularity of the mild solution.
Using  \eqref{I-spatial-temporal-S(t)} with $\mu=0, \frac{\delta_0+2}2$,  \eqref{III-spatial-temporal-S(t)} 
with $\varrho=1$ and \eqref{eq:bound-F(X)}, one can observe that, for any fixed $ \tfrac32 < \delta_0<2$,
\begin{align}\label{eq:bound-mild-solution-low-gamma}
\begin{split}
\|X(t)\|_{L^p(\Omega;\dot{H}^{\delta_0})}
\leq
&
\|E(t)X_0\|_{L^p(\Omega;\dot{H}^{\delta_0})}
\\
&
+
\int_0^t\|E(t-s) APF(X(s))\|_{L^p(\Omega;\dot{H}^{\delta_0})}\,\dd s
+
\Big(\int_0^t\|A^{\frac{\delta_0}2}E(t-r)Q^{\frac12}\|_{\mathcal{L}_2}^2\,\dd r\Big)^{\frac12}
\\
\leq&
C\Big(\|X_0\|_{L^p(\Omega;\dot{H}^{\delta_0})}
+
\int_0^t(t-s)^{-\frac{\delta_0+2}4}\|F(X(s))\|_{L^p(\Omega;H)}\,\dd s
+
\|A^{\frac{\delta_0-2}2}Q^{\frac12}\|_{\mathcal{L}_2}
\Big)
\\
\leq&
C\|X_0\|_{L^p(\Omega;\dot{H}^{\delta_0})}
+
C\sup_{s\in[0,T]}\|F(X(s))\|_{L^p(\Omega;H)}
+
C\|Q^{\frac12}\|_{\mathcal{L}_2}
<\infty,
\end{split}
\end{align}
where we also used  the Burkholder-Davis-Gundy-type inequality and the fact $Av=APv$, for any $v\in H$.
Concerning the temporal regularity of the mild solution, we  utilize    \eqref{I-spatial-temporal-S(t)}-\eqref{III-spatial-temporal-S(t)},  \eqref{eq:bound-F(X)}
and   the Burkholder-Davis-Gundy-type inequality to get, for $\beta\in[0,\delta_0]$ with  $ \tfrac32 < \delta_0<2$,
\begin{align}\label{eq:tmeporal-regualirty-H}
\begin{split}
\|X(t)-&X(s)\|_{L^p(\Omega; \dot{H}^\beta)}
\leq
\|(E(t-s)-I)X(s)\|_{L^p(\Omega; \dot{H}^\beta)}
\\
&
+
\int_s^t\|E(t-r)APF(X(r))\|_{L^p(\Omega; \dot{H}^\beta)}\,\dd r
+
C\Big(\int_s^t\|A^{\frac\beta2}E(t-r)Q^{\frac12}\|_{\mathcal{L}_2}^2\,\dd r\Big)^{\frac 12}
\\
\leq
&
C(t-s)^{\frac{\delta_0-\beta}4}\Big(\| X(s)\|_{L^p(\Omega;\dot{H}^{\delta_0})}
+
\|A^{\frac{\delta_0-2}2}Q^{\frac12}\|_{\mathcal{L}_2}\Big)
+
C\int_s^t(t-r)^{-\frac{2+\beta}4}\|F(X(r))\|_{L^p(\Omega; H)}\,\dd r
\\
\leq
&
C(t-s)^{\frac{\delta_0-\beta}4}
\Big(
\sup_{s\in[0,T]}\| X(s)\|_{L^p(\Omega; \dot{H}^{\delta_0})}
+
\|Q^{\frac12}\|_{\mathcal{L}_2}+(t-s)^{\frac{2-\delta_0}4}
\sup_{s\in[0,T]} \| F(X( s )) \|_{L^p(\Omega; H)}
\Big)
\\
\leq
&
C(t-s)^{\frac{\delta_0-\beta}4},
\end{split}
\end{align}
where we also used the fact that $ (t-s)^{\frac{2-\delta_0}4} \leq T^{\frac{2-\delta_0}4}$ for $\tfrac32 < \delta_0<2$, 
$0 \leq s<t \leq T$ and thus the final constant $C$ depends on $T$.
Next we show an improved spatial regularity of the mild solution.   First, the above two estimates and \eqref{eq:embedding-equatlity-IIII}  imply,
\begin{align}\label{eq:bound-F(x)-F(y)-H(-1)}
\begin{split}
& \|A^{-\frac12}P\big(F(X(s))-F(X(t))\big)\|_{L^p(\Omega;\dot{H})}
\leq
C
\|F(X(s))-F(X(t))\|_{L^p(\Omega;L_{\frac65})}
\\
&\quad \leq
C
\|X(s)-X(t)\|_{L^{2p}(\Omega;\dot{H})}
\Big(
1 + \sup_{s\in[0,T]}\|X(s)\|_{L^{4p}(\Omega;L_6)}^2 
\Big)
\\
&\quad \leq
C|t-s|^{\frac{\delta_0}4
},
\quad
\forall \, \delta_0\in(\tfrac32,2).
\end{split}
\end{align}
Combining this with \eqref{eq:bound-F(X)}, \eqref{I-spatial-temporal-S(t)}-\eqref{IV-spatial-temporal-S(t)}  and
 the Burkholder-Davis-Gundy-type inequality shows, for $\delta_0\in(\frac32,2)$ and $t \in [0, T]$,
\begin{align}\label{eq:H2-bound-mild-solution}
\begin{split}
\|X(t)\|_{L^p(\Omega;\dot{H}^2)}
\leq
&
\|E(t)X_0\|_{L^p(\Omega;\dot{H}^2)}
+
\int_0^t\|E(t-s)A^2P\big(\,F(X(s))-F(X(t))\,\big)\|_{L^p(\Omega;\dot{H})}\,\dd s
\\
&+
\Big\|\int_0^tE(t-s)A^2PF(X(t))\,\dd s\Big\|_{L^p(\Omega;\dot{H})}
+
\Big\|\int_0^tAE(t-s)\,\dd W(s)\Big\|_{L^p(\Omega;\dot{H})}
\\
\leq
&
C\|X_0\|_{L^p(\Omega;\dot{H}^2)}
+
C\int_0^t(t-s)^{-\frac54}\|A^{-\frac12}P\big(F(X(s))-F(X(t))\big)\|_{L^p(\Omega;\dot{H})}\,\dd s
\\
&
+
C\|F(X(t))\|_{L^p(\Omega;H)}
+
C
\left(\int_0^t\|AE(t-s)Q^{\frac12}\|_{\mathcal{L}_2}^2\,\dd s\right)^{\frac 12}
\\
\leq
&
C
\Big(
\|X_0\|_{L^p(\Omega;\dot{H}^2)}
+
\int_0^t(t-s)^{\frac{\delta_0-5}4}\,\dd s
+
\sup_{s\in[0,T]}\|F(X(s))\|_{L^p(\Omega;H)}
+
\|Q^{\frac12}\|_{\mathcal{L}_2}
\Big)
\\
<
&
\infty.
\end{split}
\end{align}
Taking the above estimate and \eqref{eq:algebra-properties-Hs} into account, one deduces
\begin{align}\label{bound-AF}
\sup_{s\in[0,T]}\| PF(X(s)) \|_{L^p(\Omega;\dot{H}^2)}
\leq
C\sup_{s\in[0,T]}\| P F(X(s)) \|_{L^p(\Omega;H^2(D))}
\leq
C
\Big(
1
+
\sup_{s\in[0,T]}\|X(s)\|^3_{L^{3p}(\Omega;\dot{H}^2)}
\Big)
<
\infty,
\end{align}
where we recalled $ P v = v - |D|^{-1}\int_D v\,\dd x, v \in H$.
Then, repeating the same lines in the proof of \eqref{eq:bound-mild-solution-low-gamma} and \eqref{eq:tmeporal-regualirty-H}
can readily result in
\begin{align*}
\sup_{s\in[0,T]}\|X(s)\|_{L^p(\Omega;\dot{H}^3)}
<\infty,
\end{align*}
and
\begin{align*}
\|X(t)-X(s)\|_{L^p(\Omega;\dot{H}^2)}
\leq
C|t-s|^{\frac14},
\quad
\forall \, 0 \leq s < t \leq T.
\end{align*}
Further, the above estimates together with \eqref{eq:algebra-properties-Hs} and \eqref{eq:embedding-equatlity-I} imply
\begin{align}
\|P(F(X(t))-F(X(s)))\|_{L^p(\Omega;\dot{H}^2)}
&\leq
C\|P ( F(X(t))-F(X(s)) )\|_{L^p(\Omega;H^2(D))}
\nonumber 
\\
&\leq
C\|X(t)-X(s)\|_{L^{2p}(\Omega;\dot{H}^2)}
\Big(
1+ \sup_{s\in[0,T]} \|X(s)\|^2_{L^{4p}(\Omega;\dot{H}^2)}
\Big)
\nonumber
\\
&
\leq
C|t-s|^{\frac14}.
\label{eq:error-Fx-H2}
\end{align}
Bearing this in mind and applying \eqref{I-spatial-temporal-S(t)} with $\mu=\frac\beta2$, 
\eqref{IV-spatial-temporal-S(t)} with $\rho=\frac\beta4$ and \eqref{bound-AF}, one can prove, for $\beta\in[0,4]$ and $0 \leq s < t \leq T$,
\begin{align}\label{eq:F-integrand}
\begin{split}
\Big\|\int_s^tE(t-r)APF(X(r))\,\dd r\Big\|_{L^p(\Omega;\dot{H}^\beta)}
\leq
&
\int_s^t\|E(t-r) A^{\frac{\beta}2}\|_{\mathcal{L}(\dot{H})}\|P(F(X(r))-F(X(t)))\|_{L^p(\Omega;\dot{H}^2)}\,\dd r
\\
&+
\Big\|\int_s^tE(t-r) A^{\frac{\beta+2}2} P F(X(t))\,\dd r \Big\|_{L^p(\Omega;\dot{H})}
\\
\leq
&
C\int_s^t(t-r)^{\frac{1-\beta}4}\,\dd r
+
C(t-s)^{\frac{4-\beta}4}\| P F ( X(t) )\|_{L^p(\Omega;\dot{H}^2)}
\\
\leq
&
C(t-s)^{\frac{4-\beta}4}.
\end{split}
\end{align}
This together with \eqref{I-spatial-temporal-S(t)}-\eqref{IV-spatial-temporal-S(t)} and the Burkholder--Davis-Gundy inequality gives
\begin{align}\label{eq:bound-mild-solution--gamma}
\begin{split}
\|X(t)\|_{L^p(\Omega;\dot{H}^\gamma)}
\leq
&
\|E(t)X_0\|_{L^p(\Omega;\dot{H}^\gamma)}
+
\Big\|\int_0^tE(t-s) APF(X(s))\,\dd s\Big\|_{L^p(\Omega;\dot{H}^\gamma)}
\\
&
+
\Big(\int_0^t\|A^{\frac\gamma2} E(t - s )Q^{\frac12}\|_{\mathcal{L}_2}^2\,\dd s \Big)^{\frac12}
\\
\leq&
C\Big(\|X_0\|_{L^p(\Omega;\dot{H}^\gamma)}
+
t^{\frac{4-\gamma}4}
+
\|A^{\frac{\gamma-2}2}Q^{\frac12}\|_{\mathcal{L}_2}
\Big)
<\infty,
\quad
\forall \, t \in [0, T].
\end{split}
\end{align}
This confirms \eqref{them:spatial-regularity-mild-solution}.
To prove  \eqref{them:temporal-regularity-mild-stoch}, we use
\eqref{eq:F-integrand}, \eqref{eq:bound-mild-solution--gamma},  \eqref{II-spatial-temporal-S(t)} with $\nu=\frac{\gamma-\beta}2$ and \eqref{III-spatial-temporal-S(t)} with $\varrho = \max\{ \frac{\beta+2-\gamma}2, 0\}$ to derive that
\begin{align}\label{eq:proof-temporal-beta-in-
gamma-2}
\begin{split}
\|X(t)-&X(s)\|_{L^p(\Omega; \dot{H}^\beta)}
\leq
\|(E(t-s)-I)X(s)\|_{L^p(\Omega; \dot{H}^\beta)}
\\
&
+
\Big\|\int_s^tE(t-r)A P F ( X(r) ) \,\dd r\Big\|_{L^p(\Omega; \dot{H}^\beta)}
+
C\Big(\int_s^t \| A^{\frac{\beta - \gamma + 2} 2} E( t - r ) A^{\frac{\gamma - 2} 2} Q^{\frac12}\|_{\mathcal{L}_2}^2\,\dd r\Big)^{\frac 12}
\\
\leq
&
C(t-s)^{\frac{\gamma-\beta}4}\| X(s)\|_{L^p(\Omega;\dot{H}^\gamma)}
+
C(t-s)^{\frac{4-\beta}4}
+
C( t - s)^{\frac12 [ 1 - \max \{ \frac{\beta+2-\gamma}2, 0 \} ] }\|A^{\frac{\gamma-2}2}Q^{\frac12}\|_{\mathcal{L}_2}
\\
\leq
&
C(t-s)^{\frac{\gamma-\beta}4}\| X(s)\|_{L^p(\Omega;\dot{H}^\gamma)}
+
CT^{\frac{4 - \gamma}{4} } (t-s)^{\frac{\gamma-\beta}4}
+
C( t - s)^{\min \{ \frac{\gamma - \beta}{4}, \frac12 \} }\|A^{\frac{\gamma-2}2}Q^{\frac12}\|_{\mathcal{L}_2}
\\
\leq
&
C(t-s)^{ \min \{ \frac{\gamma - \beta}{4}, \frac12 \} },
\end{split}
\end{align}
as required.
This thus finishes the proof of the theorem.
$\square$

As a direct consequence of Theorem \ref{them:regulairty-mild-solution}, the following lemma holds.
\begin{lemma}\label{lem:sptial-temporal-F}
Let Assumptions \ref{assum:linear-operator-A}-\ref{assum:intial-value-data} be fulfilled. Then the following results hold
\begin{align}\label{lem:bound-A-deriveaton-f}
\sup_{s\in[0,T]}
\|
A^{\frac12}
P F'( X ( s ) ) A^{ -\frac12 } P v \|_{L^p(\Omega; \dot{H} )}
\leq
C
\Big(
1
+
\sup_{s\in[0,T]}
\|
X ( s )
\|^2_{ L^{2p}(\Omega; \dot{H}^2 ) }
\Big)
\|
v
\|_{L_6},
\:
\forall p\geq 1, v \in L_6 (D),
\end{align}
\begin{align}\label{lem:error-F(X)-F(Y)}
\|P\big(F(X(t))-F(X(s))\big)\|_{L^p(\Omega;\dot{H}^\beta)}
\leq
C|t-s|^{\min\{\frac12,\frac{\gamma-\beta}4\}}
,\;
\forall p \geq 1, 0\leq s<t\leq T,
\,
\beta\in\{0,1,2\},
\end{align}
where $\gamma\in[3,4]$ comes from Assumption \ref{assum:noise-term}.
\end{lemma}

{\it Proof of Lemma \ref{lem:sptial-temporal-F}.}
Note first that $f'(v)=3 v^2 - 1, v \in \R$.  Thus, from \eqref{eq:embedding-equatlity-I}, \eqref{them:spatial-regularity-mild-solution} and H\"older's inequality, it follows that, for any $ v \in L_6 (D)$,
\begin{align*}
\begin{split}
&
\|
A^{\frac12} P F'( X(s) ) A^{-\frac12} P v
\|_{L^p(\Omega; \dot{H} )}
\leq
\|
\nabla \big( 3X^2(s) - 1 \big) A^{-\frac12} P v
\|_{L^p(\Omega;H)}
+
\|
( 3X^2(s) - 1 ) \nabla A^{-\frac12} P v
\|_{L^p(\Omega;H)}
\\
&\quad \leq
C
\big(
1
+
\|X ( s ) \nabla X(s)\|_{L^{p}(\Omega;L_3)} + \| X^2 ( s ) \|_{L^{p}(\Omega; V)}
\big)
( \| A^{-\frac12} P v \|_{L_6} + \| v \|_{L_6} )
\\
&\quad \leq
C
\big(
1
+
\|X(s)\|_{L^{2p}(\Omega;\dot{H}^2)}^2
\big)\| v \|_{L_6}
.
\end{split}
\end{align*}
To validate \eqref{lem:error-F(X)-F(Y)}, we first recall \eqref{them:temporal-regularity-mild-stoch} 
and \eqref{eq:error-Fx-H2} to attain the desired assertion for the case $\beta = 2$.
With regard to the case $\beta = 1$, one can
apply \eqref{eq:embedding-equatlity-III},  Sobolev's inequality $\|u\|_{L_3}\leq C \|u\|_{H^1 (D)}
\leq C |u|_1$  and  H\"older's inequality to show
\begin{align}\label{eq:bound-F-H1-semi}
\begin{split}
& \|A^{\frac12}P\big(F(X(t) ) - F ( X(s))\big)\|
\leq
\big\|
\nabla\big(F(X(t,\cdot))-F(X(s,\cdot))\big)
\big\|
\\
& \quad \leq
\|\big(X(t,\cdot)-X(s,\cdot)\big)\cdot\nabla\big(X^2(t,\cdot)+X(t,\cdot) X(s,\cdot)+X^2(s,\cdot)\big)\|
\\
&\qquad +
\|\nabla\big(X(t,\cdot)-X(s,\cdot)\big)\cdot \big(X^2(t,\cdot)+X(t,\cdot)X(s,\cdot)+X^2(s,\cdot)\big)\|
+
|X(t)-X(s)|_1
\\
& \quad
\leq
C\|X(t)-X(s)\|_{L_6}(\|\nabla X(t)\|_{L_3}+\|\nabla X(s)\|_{L_3}) (\|X(t)\|_V+\|X(s)\|_V)
\\
&
\qquad +
C| X(t)-X(s)|_1(1+\|X(t)\|_V^2+\|X(s)\|_V^2)
\\
&
\quad \leq
C|X(t)-X(s)|_1(1+|X(t)|_2^2+|X(s)|_2^2).
\end{split}
\end{align}
Further, combining this with  \eqref{them:temporal-regularity-mild-stoch}
enables us to obtain
\begin{align*}
\begin{split}
\|P\big(F(X(t))-F(X(s))\big)\|_{L^p(\Omega;\dot{H}^1)}
& \leq
C\|X(t)-X(s)\|_{L^{2p}(\Omega;\dot{H}^1)}
\Big(1+
\|X(t)\|^2_{L^{4p}(\Omega;\dot{H}^2)}
+
\|X(s)\|^2_{L^{4p}(\Omega;\dot{H}^2)}
\Big)
\\
&
\leq
C|t-s|^{\min\{\frac12,\frac{\gamma-1}4\}},
\end{split}
\end{align*}
verifying \eqref{lem:error-F(X)-F(Y)} for the case $\beta = 1$. Similarly,  one can easily deduce  \eqref{lem:error-F(X)-F(Y)}  for $\beta = 0$,
by taking \eqref{eq:local-condition} and \eqref{them:temporal-regularity-mild-stoch} into account.
Hence the proof of this lemma is complete.
$\square$

\section{The finite element spatial semi-discretization}
In this section, we examine the finite element spatial semi-discretization of the CHC equation and show uniform-in-time moment bounds
of the solution to the semi-discrete problem, which will be used later for the convergence analysis.
\subsection{Basic elements of the finite element spatial discretization}

Before coming to the semi-discrete finite element method (FEM) for \eqref{eq:abstract-SCHE},
we make the following assumptions.
\begin{assumption}
\label{assumption:triangulations}
Suppose that the spatial domain $D \subset \mathbb{R}^d$, $d\in\{1,2,3\}$ has polygonal boundary $\partial D$.
The triangulations $\{ \mathcal{T}_h\}_{h >0}$ of $D$ with maximal mesh size $h$, are assumed to be quasi-uniform.
Let $\{V_h\}_{h>0}\subset H^1(D)$ be the space of continuous functions that are piecewise polynomials of degree at most $r-1$
over $\mathcal{T}_h$, for $r\in\{2,3,4\}$ in dimension $d=1$ and $r=2$ in dimension $d\in\{2,3\}$.
\end{assumption}

Furthermore, 
we  define $\dot{V}_h=PV_h$ by
\[
\dot{V}_h=\Big\{v_h\in V_h: \int_Dv_h\,\dd x=0\Big\},
\]
and introduce a discrete Laplace operator $A_h:V_h\rightarrow V_h$ defined by
\begin{align}\label{eq:definition-discrete-A}
(A_hv_h,\chi_h)=a(v_h,\chi_h):=(\nabla v_h, \nabla \chi_h),\quad \forall v_h\in V_h,\;\chi_h\in V_h.
\end{align}
The operator $A_h$ is selfadjoint, positive semidefinite on
 $V_h$ and positive definite on $\dot{V}_h$, and has an orthonormal eigenbasis $\{e_{j,h}\}_{j=0}^{\mathcal{N}_h}$ in ${V}_h$ with
 corresponding eigenvalues $\big\{\lambda_{j,h}\big\}_{j=0}^{\mathcal{N}_h}$, satisfying
\begin{align*}
0 =  \lambda_{0,h} < \lambda_{1,h} \leq \cdots \leq \lambda_{j,h}\leq\cdots\leq \lambda_{\mathcal{N}_h, h},
\end{align*}
where $\mathcal{N}_h := \dim(V_h)$ and $e_{0,h} = e_0 = |D|^{-\frac12}$.
Also, we introduce a discrete norm  on $\dot{V}_h$, defined by
\begin{align*}
|v_h|_{\alpha,h}=\|A_h^{\frac\alpha2}v_h\|
=
\Big(\sum_{j=1}^{\mathcal{N}_h}\lambda_{j,h}^\alpha|
(v_h,e_{j,h})|^2\Big)^{\frac12}, \;v_h\in \dot{V}_h,\;\alpha\in \mathbb{R},
\end{align*}
which corresponds to the discrete inner product $(v,w)_{\alpha,h}:=(A_h^{\frac{\alpha}{2}} v,  A_h^{\frac{\alpha}{2}}  w)$, 
$\forall v,w\in \dot{V}_h$. Note that
\begin{align}\label{eq-eqvilent-discrete-norem-H1-norm}
|v_h|_1=\|A^{\frac12}v_h\|=\|\nabla v_h\|=\|A_h^{\frac12}v_h\|=|v_h|_{1,h},
\quad v_h\in\dot{ V}_h.
\end{align}
In addition, we introduce a Riesz representation operator $R_h:\dot{H}^1\rightarrow \dot{V}_h$ defined by
\begin{align}\label{eq:definition-Rh}
a(R_hv,\chi_h)=a(v,\chi_h),
\quad
\forall v\in \dot{H}^1,
\,
\chi_h\in \dot{V}_h,
\end{align}
and a generalized projection operator $P_h:H\rightarrow V_h$ given by
\begin{align}\label{eq:definition-Ph}
(P_h v,\,\chi_h)=(v,\,\chi_h),
\quad
\forall v\in H,
\,
\chi_h\in V_h .
\end{align}
It is easy to see that $P_h$  is also a projection operator from $\dot{H}$ to $\dot{V}_h$ and
\begin{align}\label{eq:relation-A-Ah-Rh-Ph}
 P_hA=A_hR_h.
 \end{align}
Owing to Assumption \ref{assumption:triangulations}, we have the following error bounds for the operators $R_h$ and $P_h$ 
(cf. \cite[Chapter 1]{ee2006galerkin} and \cite[(2.3)]{larsson2011finite}),
\begin{align}\label{eq:error-interpolation}
\begin{split}
|(I-R_h) v|_i+|(I-P_h) v|_i
\leq
 Ch^{\beta-i}|v|_\beta,
 \;
 \forall v\in \dot{H}^\beta,\;i=0,1,\;\beta\in[1,r],
\end{split}
\end{align}
where $r\in\{2,3,4\}$ for $d=1$ and $r=2$ for $d\in\{2,3\}$.
We mention that only $r = 2$ is considered for $d\in\{2,3\}$, because for higher-order elements, 
i.e., $ r \in \{ 3, 4\}$, the situation becomes very complicated in high dimension $d \in \{ 2, 3\}$.
%
%
Also, the quasi-uniform mesh $\mathcal{T}_h$  ensures that $P_h$ is bounded with respect to the $\dot{H}^1$ and $L_4$ norms and that we have an inverse bound for $A_h$,
\begin{align}
\|P_hv\|_{L_4}
&\leq
C\|v\|_{L_4},\;\forall v\in L_4,
\label{eq:bound-Ph-L4}
\\
|P_hv|_1
&\leq
 C|v|_1, \;\;\forall v\in \dot{H}^1,
 \label{eq:bound-Ph-H1}
 \\
\|A_h v_h \|
&\leq
C h^{-2} \| v_h \|, \;\;\forall v_h \in V_h.
\label{eq:inverse-ineq}
\end{align}
These three inequalities have been also claimed in \cite[(2.12)]{kovacs2011finite}.
The inverse inequality \eqref{eq:inverse-ineq} together with \eqref{eq:error-interpolation} helps us to infer
\begin{align}
\label{eq:bound-Ph-H2}
\|A_hP_hv\|
\leq
\|A_hP_h(I-R_h)v\|
+
\|P_hAv\|
\leq
Ch^{-2}\|(I-R_h)v\|
+
C|v|_2
\leq
C|v|_2,
\quad
\forall v \in \dot{H}^2.
\end{align}
Moreover, the operators $A$ and $A_h$ obey
\begin{align}\label{eq:relation-A-Ah}
C_1\|A_h^{\frac \alpha2}P_hv\|
\leq
\|A^{\frac \alpha2}v\|
\leq
C_2
\|A_h^{\frac \alpha2}P_hv\|,
\quad
\forall
v\in \dot{H}^{\alpha},
\:
\alpha\in[-1,1],
\end{align}
and similarly to that  in \cite[Theorem 6.11, (6.91)]{ee2006galerkin},
\begin{align}\label{eq:V-norm-control-by-Ah-norm}
\|v_h\|_V\leq C\|A_h v_h\|,
\quad
\forall v_h\in \dot{V}_h.
\end{align}
Combining \eqref{eq:bound-Ph-H2} with \eqref{eq:relation-A-Ah} gives
\begin{equation}
\label{eq:Ah-A-bound}
\|A_h^{\frac \alpha2}P_hv\|
\leq
C \|A^{\frac \alpha2}v\|,
\quad
\forall
v\in \dot{H}^{\alpha},
\:
\alpha \in [-1, 2].
\end{equation}

\subsection{Moment bounds of the approximation}
In this subsection,  we come to the semi-discrete finite element approximation of the stochastic
 problem and provide some useful moment bounds for the semi-discrete approximations.

The semi-discrete finite element method for the problem \eqref{eq:abstract-SCHE} can be written as,
\begin{align}\label{eq:semi-discrete-problem}
\,\dd X_h(t)+A_h\big(A_hX_h(t)+ P_hF(X_h(t))\big)\,\dd t=P_h\,\dd W(t),
\quad
X_h(0)=P_hX_0.
\end{align}
Similarly to \eqref{eq:definition-semigroup-E(t)}, the analytic semigroup  $E_h(t)$ generated by the discrete operator $-A^2_h$ can be given as follows,
\begin{align*}
E_h(t)P_hv=e^{-tA_h^2}P_hv=\sum_{j=0}^{ \mathcal{N}_h } e^{-t\lambda_{j,h}^2}(P_hv,e_{j,h})e_{j,h}=PE_h(t)P_hv+(I-P)v.
\end{align*}
Since $\dot{V}_h$ is finite-dimensional and $F$ is a polynomial of particular structure and recalling $X_0, W(t)$ are $\dot{H}$-valued and thus $P_hX_0, P_h W (t)$ are $\dot{V}_h$-valued, one can easily check that the  problem \eqref{eq:semi-discrete-problem} admits a unique $\dot{V}_h$-valued solution $X_h(t)$, adapted, almost surely continuous, satisfying
\begin{align*}
X_h(t)-P_hX_0
+
\int_0^t\big(A_h^2X_h(s)+A_hP_hF(X_h(s))\big)\,\dd s
=
P_hW(t),
\end{align*}
or equivalently in a mild solution form,
\begin{align}\label{eq:mild-solution-semi-discrete-problem}
X_h(t)=E_h(t)P_hX_0
-
\int_0^tE_h(t-s)A_hP_hF(X_h(s))\,\dd s
+
\int_0^tE_h(t-s)P_h\,\dd W(s).
\end{align}
%
%
Before presenting moment bounds of the approximations,
we introduce a spatially discrete version of \eqref{I-spatial-temporal-S(t)}-\eqref{IV-spatial-temporal-S(t)},
 which plays an important role in deriving the moment bounds of $X_h$.
\begin{lemma}\label{lem:propertiy-semgroup-semi}
Under Assumptions \ref{assum:linear-operator-A}, \ref{assumption:triangulations}, the following estimates for $E_h$ hold,
\begin{align}\label{lem:eq-smooth-property-Eh(t)}
\|A_h^{\mu}E_h(t)P_hv\|
&\leq
Ct^{-\frac\mu2}\|v\|,
\quad
\forall \mu\geq 0, \: v\in \dot{H},
\\
\label{lem:eq-temporal-smooth-property-Eh(t)}
\|A_h^{-\nu}(I-E_h(t))P_hv\|
 &\leq
 Ct^{\frac\nu2}\|v\|,
 \quad
 \forall \nu\in[0,2], \: v\in \dot{H},
\\
\label{lem:eq-integrand-smooth-property-Eh(t)}
\Big\|\int_0^tA_h^2E_h(s)P_hv\,\mathrm{d}s \Big\|
& \leq
C\|v\|,
\quad
\forall v\in \dot{H},
\\
\label{lem:eq-integrand-smooth-property-Eh(t)-II}
\Big(\int_0^t\|A_h E_h(s)P_hv\|^2\,\mathrm{d}s \Big)^{\frac12}
&\leq
C\|v\|,
\quad
\forall v\in \dot{H}.
\end{align}
\end{lemma}
Note that estimates \eqref{lem:eq-smooth-property-Eh(t)}, \eqref{lem:eq-integrand-smooth-property-Eh(t)-II} 
can be found in \cite[(2.1)-(2.2)]{larsson2011finite}.
Actually, thanks to the expansion of the eigenbasis of $A_h$  in $\dot{V}_h$ and  the Parseval identity,
one can follow standard arguments to derive these estimates.
Next, we are ready to show the following moment bounds for the FEM approximation.
\begin{theorem}\label{lem:bound-solution-semi}
Let $X_h(t)$ be the  solution to \eqref{eq:semi-discrete-problem}.
If Assumptions \ref{assum:linear-operator-A}-\ref{assum:intial-value-data}, \ref{assumption:triangulations} are valid,
then 
\begin{align}\label{lem:eq-bound-solution-semi}
\sup_{s\in[0,T]}\|A_h X_h( s ) \|_{L^p(\Omega;\dot{H})}
+
\Big\|\int_0^T|A_hX_h(s)+P_hF(X_h(s))|_1^2\,\mathrm{d}s\Big\|_{L^p(\Omega;\mathbb{R})}
<\infty,
\quad
\forall p\geq 1.
\end{align}
\end{theorem}

{\it Proof of Theorem \ref{lem:bound-solution-semi}.}
A slight modification of the proof of \cite[Theorem 3.1]{kovacs2011finite} enables us to obtain
\begin{align*}
\begin{split}
&
\Big\|
\sup_{s\in[0,T]} J(X_h(s))
\Big\|^p_{L^p(\Omega;\mathbb{R})}
+
\Big\|\int_0^T|A_hX_h(s)+P_hF(X_h(s))|_1^2\,\mathrm{d}s\Big\|^p_{L^p(\Omega;\mathbb{R})}
\\
& \quad \leq
C
\bigg(
1+\|J(P_hX_0)\|^p_{L^p(\Omega;\mathbb{R})}
+
\Big\|
\sup_{t\in[0,T]}\Big|\int_0^t\big(A_hX_h(s)+P_hF(X_h(s)), P_h\,\dd W(s)\,\big)\Big|^p\Big\|_{L^1(\Omega;\mathbb{R})}
\bigg),
\end{split}
\end{align*}
where $J$ is defined by \eqref{eq:definition-lyapunov-fun}.
Owing to \eqref{eq:bound-Ph-L4} and \eqref{eq:bound-Ph-H1}, one knows $\|J(P_hX_0)\|_{L^p(\Omega;\mathbb{R})}<\infty$. 
Further, with the aid of  the Burkhoder-Davis-Gundy-type inequality, one can find that
\begin{align}\label{eq:J-sup-AX+PhF}
\begin{split}
&
\Big\|
\sup_{s\in[0,T]} J(X_h(s))
\Big\|^p_{L^p(\Omega;\mathbb{R})}
+
\Big\|\int_0^T|A_hX_h(s)+P_hF(X_h(s))|_1^2\,\mathrm{d}s\Big\|^p_{L^p(\Omega;\mathbb{R})}
\\
&\quad \leq
C\big(1
+
\Big\|\int_0^T\|Q^{\frac12}\big(A_hX_h(s)+P_hF(X_h(s))\,\big)\|^2\,\dd s\Big\|^{\frac p2}_{L^{\frac p2}(\Omega;\mathbb{R})}
\big)
\\
&\quad \leq
C\big(1
+
\|Q^{\frac12}\|_{\mathcal{L}(\dot{H})}^p
\Big\|\int_0^T\|A_hX_h(s)+P_hF(X_h(s))\|^2\,\dd s\Big\|^{\frac p2}_{L^{\frac p2}(\Omega;\mathbb{R})}\big)
\\
& \quad
\leq
C\Big(1
+
\tfrac {\|Q^{\frac 12}\|_{\mathcal{L}(\dot{H})}^{2p}}{2\varepsilon}+\frac {\varepsilon}2
\Big\|\int_0^T|A_hX_h(s)+P_hF(X_h(s))|_1^2\,\dd s\Big\|^p_{L^p(\Omega;\mathbb{R})}\Big),
\end{split}
\end{align}
where we also used  the fact
$\|Q^{\frac12}\|_{\mathcal{L}(\dot{H})}\leq \|Q^{\frac12}\|_{\mathcal{L}_2}<\infty$.
Taking  $\varepsilon>0$ small enough in \eqref{eq:J-sup-AX+PhF}, we conclude that
\begin{align*}
\Big \|
\sup_{s\in[0,T]} J(X_h(s))
\Big\|^p_{L^p(\Omega;\mathbb{R})}
+
\Big\|\int_0^T|A_hX_h(s)+P_hF(X_h(s))|_1^2\,\mathrm{d}s\Big\|^p_{L^p(\Omega;\mathbb{R})}<\infty.
\end{align*}
It remains to bound $\|A_hX_h(t)\|_{L^p(\Omega;\dot{H})}$. From the definition of the Lyapunov functional $J(\cdot)$ and noting $\Phi(s)=\frac14(s^2-1)^2$,   one can deduce that
\begin{align}\label{eq:relation-H1-J}
|v|_1^2
\leq
2J(v),
\quad
\forall v\in \dot{H}^1,
\end{align}
which leads to
\begin{align}\label{eq:bound-semidiscrete-mild-H1}
\sup_{s\in[0,T]} \|X_h(s)\|_{L^p(\Omega;\dot{H}^1)}
\leq
C\big\|\sup_{s\in[0,T]} \big[J(X_h(s))\big]^{\frac p2}\big\|_{L^1(\Omega;\mathbb{R})}^{\frac1p}< \infty.
\end{align}
Using \eqref{eq:embedding-equatlity-I} shows
\begin{align}\label{eq:bound-semi-F}
\sup_{s\in[0,T]}\|F(X_h(s))\|_{L^p(\Omega; H)}
\leq
C\Big(\sup_{s\in[0,T]}\|X_h(s)\|_{L^p(\Omega; \dot{H})}
+
\sup_{s\in [0,T]}\|X_h(s)\|_{L^{3p}(\Omega; \dot{H}^1)}^3 \Big)
<
\infty.
\end{align}
With the above estimate, one can follow the same lines of the proof of \eqref{eq:bound-mild-solution-low-gamma} to show,
for $\delta_0\in(\frac3 2, 2)$,
\begin{align*}
\begin{split}
\sup_{s\in[0,T]}
\big\|
A_h^{\frac{\delta_0}2}X_h(s)
\big\|_{L^p(\Omega; \dot{H})}
&\leq
C
\Big(
\|A_h^{\frac{\delta_0}2}P_hX_0\|_{L^p(\Omega; \dot{H})}
+
\sup_{s\in[0,T]}\|F(X_h(s))\|_{L^p(\Omega; H)}
+
\|Q^{\frac12}\|_{\mathcal{L}_2}
\Big)
\\
&\leq
C
\big(
\|X_0\|_{L^p(\Omega; \dot{H}^2)}+1
\big)
<\infty,
\end{split}
\end{align*}
where we also used \eqref{eq:bound-Ph-H2} in the second step.
Similarly to \eqref{eq:tmeporal-regualirty-H} in the previous proof, one can show
\begin{align}\label{eq:temporal-full-solution}
\|X_h(t)-X_h(s)\|_{L^p(\Omega;\dot{H})}
\leq
C |t-s|^{\frac{\delta_0}4},
\quad
\delta_0 \in ( \tfrac32, 2 ),
\end{align}
which together with \eqref{eq:embedding-equatlity-IIII} and \eqref{eq:relation-A-Ah} yields,
\begin{align}\label{eq:regularity-F-full-H(-1)}
\begin{split}
\|A_h^{-\frac12}P_hP\big(F(X_h(s))-F(X_h(t))\big)\|_{L^p(\Omega;\dot{H})}
&\leq
C\|A^{-\frac12}P\big(F(X_h(s))-F(X_h(t))\big)\|_{L^p(\Omega;\dot{H})}
\\
&\leq
C\|F(X_h(s))-F(X_h(t))\|_{L^p(\Omega;L_{\frac65})}
\\
&
\leq
C\|X_h(s)-X_h(t)\|_{L^{2p}(\Omega;\dot{H})}
\Big ( 1 +   \sup_{s\in[0,T]} \| X_h(s) \|_{L^{4p}(\Omega;L_6)}^2 \,\Big)
\\
&
\leq
C|t-s|^{\frac{\delta_0}4},
\quad
\delta_0 \in ( \tfrac32, 2 ).
\end{split}
\end{align}
Combining this  with \eqref{eq:bound-Ph-H2}, \eqref{eq:bound-semi-F}, \eqref{lem:eq-integrand-smooth-property-Eh(t)}, \eqref{lem:eq-integrand-smooth-property-Eh(t)-II} and \eqref{lem:eq-smooth-property-Eh(t)} with $\mu=0, \frac52$ gives,  for $\delta_0\in(\frac32,2)$,
\begin{align}\label{eq:H2-bound-mild-solution-smei}
\begin{split}
\|A_hX_h(t)\|_{L^p(\Omega;\dot{H})}
\leq
&
\|A_hE_h(t)P_hX_0\|_{L^p(\Omega;\dot{H})}
+
\int_0^t\|E_h(t-s)A_h^2P_h \big(\,F(X_h(s))-F(X_h(t))\,\big)\|_{L^p(\Omega;\dot{H})}\,\dd s
\\
&+
\Big\|\int_0^tE_h(t-s)A_h^2P_h F(X_h(t))\,\dd s\Big\|_{L^p(\Omega;\dot{H})}
+
C\Big\|\int_0^tA_hE_h(t-s)P_h\,\dd W(s)\Big\|_{L^p(\Omega;\dot{H})}
\\
\leq
&
C\|A_hP_hX_0\|_{L^p(\Omega;\dot{H})}
+
C\int_0^t(t-s)^{-\frac54}\|A_h^{-\frac12}P_h \big(F(X_h(s))-F(X_h(t))\big)\|_{L^p(\Omega;\dot{H})}\,\dd s
\\
&
+
C\|F(X_h(t))\|_{L^p(\Omega;H)}
+
C\left(\int_0^t\|A_hE_h(t-s)P_hQ^{\frac12}\|_{\mathcal{L}_2}^2\,\dd s\right)^{\frac 12}
\\
\leq
&
C\big(\|X_0\|_{L^p(\Omega;\dot{H}^2)}
+
\int_0^t(t-s)^{\frac{\delta_0-5}4}\,\dd s
+
\sup_{s\in[0,T]}\|F(X_h( s ) )\|_{L^p(\Omega;H)}
+
\|Q^{\frac12}\|_{\mathcal{L}_2}\big)
<
\infty,
\end{split}
\end{align}
as required, where the Burkholder-Davis-Gundy inequality was also used.
 $\square$
\section{Strong convergence rates of the FEM semi-discretization}
The target of this part  is to derive error estimates for the semi-discrete finite element
 approximation of the stochastic problem \eqref{eq:abstract-SCHE}.
 The convergence analysis heavily relies on the moment bounds obtained in the previous section
 and error estimates for the corresponding deterministic error operators as shown below.

Define two error operators $\Psi_h $ and $\Phi_h $ for the semi-discrete approximation as follows,
\begin{align}\label{eq:definition-Psih-Phih(t)}
\Psi_h(t):=E(t)-E_h(t)P_h\quad
\text{ and }
\quad \Phi_h(t):= AE(t)-A_hE_h(t)P_h, \; t\in[0,T].
\end{align}
It is easy to check that
\begin{equation}
\label{eq:Psi-h-Phi-h-identity}
\Psi_h(t) v = \Psi_h(t) P v,
\quad
\Phi_h(t) v = \Phi_h(t) P v,
\quad
\forall v \in H,
\end{equation}
since the constant eigenmodes are cancelled.
We present in the following lemma some deterministic semi-discrete error estimates for the above two error operators.
%
%
\begin{lemma}\label{lem:deterministic-error-semi}
Under Assumptions \ref{assum:linear-operator-A}, \ref{assumption:triangulations},  the following estimates for 
$\Psi_h$ and $\Phi_h$ hold,
\begin{align}
\label{lem:eq-deterministic-error-smooth-semi}
\|\Psi_h(t)v\|
&\leq
Ch^\beta  |v|_\beta,
\quad
\forall v\in \dot{H}^\beta,
\:
\beta\in[1,r],
\\
\|\Phi_h(t)v\|
&\leq
Ch^\alpha t^{-1} |v|_{\alpha-2},
\quad
\forall v\in \dot{H}^{\alpha-2},\:
\alpha\in[1,r],
\label{lem:eq-deterministic-error-displiace-semi}
\\
\label{lem:eq-Fh-integrand}
\Big(\int_0^t\|\Psi_h(s)v\|^2\,\mathrm{d} s\Big)^{\frac12}
& \leq
Ch^{\nu} |\ln h| |v|_{\nu-2},\,
\quad
\forall v\in \dot{H}^{\nu-2},
\:
\nu\in[1,r],
\\
\label{lem:eq-chemical-potial-integrand}
\Big(\int_0^t\|\Phi_h(s)v\|^2\,\mathrm{d} s\Big)^{\frac12}
& \leq
Ch^{\mu} |\ln h| |v|_\mu,\,
\quad
\forall v\in \dot{H}^\mu, \:
\mu\in[0,r],
%
\\
\label{lem:eq-chemical-potial-integrand-II}
\Big\|\int_0^t\Phi_h(s)v\,\mathrm{d} s\Big\|
& \leq
Ch^\varrho |v|_{\varrho-2},
\quad
\forall v\in \dot{H}^{\varrho-2}, \:
\varrho\in[1,r].
\end{align}
\end{lemma}

{\it Proof of Lemma \ref{lem:deterministic-error-semi}.}
 The  estimates \eqref{lem:eq-deterministic-error-smooth-semi} and \eqref{lem:eq-Fh-integrand} are shown in
\cite[Theorem 2.1]{larsson2011finite}.
Taking \eqref{eq:error-interpolation} into account, we can make a slight modification of the proof of \cite[(5.6)]{furihata2018strong} in the case $\delta=0$ to prove \eqref{lem:eq-deterministic-error-displiace-semi}.
In order to show \eqref{lem:eq-chemical-potial-integrand}, we rely on a simple interpolation 
between the cases $\mu=0$ and $\mu=r$.
The case $\mu=0$ immediately follows
from \eqref{III-spatial-temporal-S(t)} with $\varrho=1$ and  \eqref{lem:eq-integrand-smooth-property-Eh(t)-II}.
 For the case $\mu=r$, we use \eqref{lem:eq-integrand-smooth-property-Eh(t)-II}, \eqref{eq:relation-A-Ah-Rh-Ph}, \eqref{eq:definition-Psih-Phih(t)}, \eqref{eq:error-interpolation} with $i=0$, $\beta=r$ and  \eqref{lem:eq-Fh-integrand} with $\nu=r$  to deduce
\begin{align*}
\begin{split}
\Big(\int_0^t\|\Phi_h(s)v\|^2\,\mathrm{d} s\Big)^{\frac12}
\leq
&
\Big(\int_0^t\|(E(s)-E_h(s)P_h)Av\|^2\,\mathrm{d} s\Big)^{\frac12}
+
\Big(\int_0^t\|A_hE_h(s)P_h(I-R_h)v\|^2\,\mathrm{d} s\Big)^{\frac12}
\\
\leq
&
Ch^r|\ln h| |v|_r
+
C\|(R_h-I)v\|
\leq Ch^r|\ln h||v|_r.
\end{split}
\end{align*}
Finally, an interpolation argument concludes the proof of \eqref{lem:eq-chemical-potial-integrand}.
Similarly as before, we use \eqref{eq:relation-A-Ah-Rh-Ph} to split the term $\big\|\int_0^t\Phi_h(s)v\,\mathrm{d} s\big\|$  into two parts:
\begin{align}\label{eq:Psi(s)A{-1}v-decompoese}
\begin{split}
\Big\|\int_0^t\Phi_h(s)v\,\mathrm{d} s\Big\|
&=
\Big\|\int_0^t(A^2E(s)A^{-1}-A_h^2E_h(s)R_hA^{-1})v\,\mathrm{d} s\Big\|
\\
&
\leq
\Big\|\int_0^t(A^2E(s)-A_h^2E_h(s)P_h)A^{-1}v\,\mathrm{d} s\Big\|
+
\Big\|\int_0^tA_h^2E_h(s)P_h(R_h-I)A^{-1}v\,\mathrm{d} s\Big\|
\\
&\leq
\Big\|\int_0^t\Psi'_h(s)A^{-1}v\,\mathrm{d} s\Big\|
+
\Big\|\int_0^tE'_h(s)P_h(R_h-I)A^{-1}v\,\mathrm{d} s\Big\|.
\end{split}
\end{align}
For the first term, we use the fundamental theorem of calculus,  \eqref{eq:error-interpolation} with $i=0$, $\beta=\varrho$ and \eqref{lem:eq-deterministic-error-smooth-semi} with $\beta=\varrho$ to show,  for $\varrho\in[1,r]$,
\begin{align}\label{eq:Psi(s)A{-1}v-I}
\begin{split}
\Big\|\int_0^t\Psi'_h(s)A^{-1}v\,\mathrm{d} s\Big\|
=
\|\big(\Psi_h(t)-\Psi_h(0)\big)A^{-1}v\|
\leq
\|\Psi_h(t)A^{-1}v\|+\|(I-P_h)A^{-1}v\|
\leq
Ch^\varrho |v|_{\varrho-2}.
\end{split}
\end{align}
Similarly, we combine the boundness of $ E_h(s)P_h $ in $ \dot{H} $ with \eqref{eq:error-interpolation} to yield
\begin{align}\label{eq:Psi(s)A{-1}v-II}
\begin{split}
\Big\|\int_0^tE'_h(s)P_h(R_h-I)A^{-1}v\,\mathrm{d} s\Big\|
=
\|(E_h(t)-I)P_h(R_h-I)A^{-1}v\|
\leq
Ch^\varrho |v|_{\varrho-2},
\end{split}
\end{align}
which together with \eqref{eq:Psi(s)A{-1}v-I} and \eqref{eq:Psi(s)A{-1}v-decompoese} implies \eqref{lem:eq-chemical-potial-integrand-II}. 
This finishes the proof of this lemma.
$\square$
%

At the moment, we are well-prepared to prove the main result of this section.
\begin{theorem}\label{them:error-estimates-semi-problem}
Let $X(t)$ be the weak solution of \eqref{eq:abstract-SCHE} and let $X_h(t)$  be the solution of
 \eqref{eq:semi-discrete-problem}. Also, let Assumptions \ref{assum:linear-operator-A}-\ref{assum:intial-value-data}
be valid for some $\gamma\in[3,4]$ and let Assumption \ref{assumption:triangulations} be fulfilled with 
$r\in\{2,3,4\}$ for $d=1$ and $r=2$ for $d\in\{2,3\}$.  Then for any $t \in [0, T]$ and $ p\in [1,\infty)$ it holds,
\begin{align}\label{them:eq-error-estimates-semi-problem}
\|X(t)-X_h(t)\|_{L^p(\Omega;\dot{H})}
\leq
Ch^\kappa |\ln h|,\quad \text{ with } \quad \kappa=\min\{\gamma,r\} .
\end{align}
Moreover, the discrepancy between the "chemical potential" $Y(t) : = AX(t)+PF(X(t))$ and its approximation 
$Y_h(t) :=  A_h X_h(t)+P_hPF(X_h(t))$ is measured as follows, for any $t \in (0, T]$ and $ p\in [1,\infty)$,
\begin{align}\label{them:eq-chemical potential-error-estimates-semi-problem}
\|Y(t)-Y_h(t)\|_{L^p(\Omega;\dot{H})}
\leq
C\big(1+t^{-1}\big)h^\iota |\ln h|,\quad \text{ with } \quad \iota=\min\{\gamma-2,r-1\}.
\end{align}
\end{theorem}
%
%

{\it Proof of Theorem \ref{them:error-estimates-semi-problem}.}
 Since $A$ does not commute with $P_h$, the usual arguments splitting the error $X(t)-X_h(t)$ into
$(I-P_h)X(t)$ and $P_hX(t)-X_h(t)$ do not work here. To prove this theorem,
we propose a different approach and
introduce a new auxiliary problem:
\begin{align}\label{eq:semi-discrete-auxiliary-problem}
\,\dd \widetilde{X}_h(t)+A_h\big(A_h\widetilde{X}_h(t)+P_hF(X(t))\big)\,\dd t=P_h\,\dd W(t),\, X_h(0)=P_hX_0,
\end{align}
whose unique solution can be written as, in the mild form,
\begin{align}\label{eq:mild-soluiton-auxiliary-semi}
\widetilde{X}_h(t)=E_h(t)P_hX_0
-
\int_0^tE_h(t-s)A_hP_h F(X(s))\,\dd s
+
\int_0^tE_h(t-s)P_h\,\dd W(s).
\end{align}
Now, we separate the considered error term $\|X(t)-X_h(t)\|_{L^p(\Omega;\dot{H})}$ as follows,
\begin{align}\label{eq:decompose-error-semi}
\|X(t)-X_h(t)\|_{L^p(\Omega;\dot{H})}
\leq
\|X(t)-\widetilde{X}_h(t)\|_{L^p(\Omega;\dot{H})}
+
\|\widetilde{X}_h(t)-X_h(t)\|_{L^p(\Omega;\dot{H})}.
\end{align}
Recall that a similar error decomposition was done in \cite[(5.18)]{kovacs2011finite}.
The first error term can be treated in a standard way.
Subtracting \eqref{eq:mild-soluiton-auxiliary-semi} from \eqref{them:eq-mild-solution-stochastic-equation}  yields
\begin{align}\label{eq:semi-auxiliary-error}
\begin{split}
\|X(t)-\widetilde{X}_h(t)\|_{L^p(\Omega;H)}
\leq&
\|\Psi_h(t)X_0\|_{L^p(\Omega;H)}
+
\Big\|\int_0^t\Phi_h(t-s) F(X(s))\,\dd s\Big\|_{L^p(\Omega;H)}
\\
&+
\Big\|\int_0^t\Psi_h(t-s)\,\dd W(s)\Big\|_{L^p(\Omega;H)}
\\
=:
&
I_1+I_2+I_3,
\end{split}
\end{align}
where the two error operators $\Psi_h $ and $\Phi_h $ are defined by \eqref{eq:definition-Psih-Phih(t)}.
In what follows we treat $I_1,I_2$ and $I_3$, separately. At first, we utilize \eqref{lem:eq-deterministic-error-smooth-semi} with
$\beta=\kappa$ to derive
\begin{align*}
I_1
\leq
Ch^\kappa\|X_0\|_{L^p(\Omega;\dot{H}^\kappa)},
\quad
\kappa=\min\{\gamma,r\}.
\end{align*}
Similarly, employing \eqref{bound-AF}, \eqref{lem:error-F(X)-F(Y)}, \eqref{eq:Psi-h-Phi-h-identity},
\eqref{lem:eq-deterministic-error-displiace-semi} with $\alpha=\kappa$ and \eqref{lem:eq-chemical-potial-integrand-II} with $\varrho=\kappa$
yields
\begin{align*}
\begin{split}
I_2
\leq
&
\Big \|\int_0^t\Phi_h(t-s)PF(X(t))\,\dd s\Big\|_{L^p(\Omega;\dot{H})}
 +
  \int_0^t\|\Phi_h(t-s)P\big(F(X(t))-F(X(s))\big)\|_{L^p(\Omega;\dot{H})}\,\dd s
  \\
\leq
&
Ch^{\kappa} \|PF(X(t))\|_{L^p(\Omega;\dot{H}^{ \kappa-2})}
+
Ch^{\kappa} \int_0^t(t-s)^{-1}\|P\big(F(X(t))-F(X(s))\big)\|_{L^p(\Omega;\dot{H}^{\kappa-2})}\,\dd s
\\
\leq
&
Ch^{\kappa}  \sup_{t\in[0,T]}\|PF(X(t))\|_{L^p(\Omega;\dot{H}^2)}
+
Ch^{\kappa} \int_0^t(t-s)^{-\frac34}\,\dd s
\\
\leq
&
Ch^{ \kappa},
\qquad
\kappa = \min\{\gamma,r\}.
\end{split}
\end{align*}
Finally,  we use the Burkholder-Davis-Gundy-type inequality and
\eqref{lem:eq-Fh-integrand} with $\nu=\kappa$ to arrive at
\begin{align}\label{eq:estimate-I3}
I_3\leq
C_p\Big(\int_0^t\|\Psi_h(t-s)Q^{\frac12}\|_{\mathcal{L}_2}^2\,\dd s\Big)^{\frac12}
\leq
Ch^\kappa|\ln h| \|A^{\frac{\kappa-2}2}Q^{\frac12}\|_{\mathcal{L}_2}
\leq
Ch^\kappa|\ln h|,
\quad
\kappa = \min\{ \gamma, r \}.
\end{align}
Putting the above estimates together yields
\begin{align}\label{eq:error-auxilary-semi}
\|X(t)-\widetilde{X}_h(t)\|_{L^p(\Omega;\dot{H})}\leq Ch^\kappa|\ln h|,
\qquad
\kappa=\min\{\gamma,r\}.
\end{align}
Next we turn our attention to the error $\widetilde{e}_h(t) := \widetilde{X}_h(t)-X_h(t)$, which satisfies
\begin{align}\label{eq:error-equation-semi}
\dd \widetilde{e}_h(t)
+
A_h^2\widetilde{e}_h(t) \, \dd t
=
A_h P_h \big(F(X(t)) - F(X_h(t))\,\big) \, \dd t, \quad\widetilde{e}_h(0)=0.
\end{align}
Multiplying both sides of \eqref{eq:error-equation-semi} by $A_h^{-1}\widetilde{e}_h$, 
using \eqref{eq:local-condition}, \eqref{eq:V-norm-control-by-Ah-norm}, \eqref{eq:embedding-equatlity-I} 
and recalling the fact $\|\widetilde{e}_h\|^2\leq |\widetilde{e}_h|_1|\widetilde{e}_h|_{-1,h}$ one obtains
\begin{align}\label{eq:error-widetild-e-semi}
\begin{split}
&
\tfrac12\tfrac{\dd}{\dd s}
|\widetilde{e}_h(s)|_{-1,h}^2+|\widetilde{e}_h(s)|_1^2
=
\big(F(\widetilde{X}_h(s))-F(X_h(s)),\widetilde{e}_h(s)\big)
+
\big(F(X(s))-F(\widetilde{X}_h(s)),\widetilde{e}_h(s)\big)
\\
&\quad \leq
\tfrac32\|\widetilde{e}_h(s)\|^2
+
\tfrac12\|F(X(s))-F(\widetilde{X}_h(s))\|^2
\\
&\quad \leq
\tfrac32|\widetilde{e}_h(s)|_1|\widetilde{e}_h(s)|_{-1,h}
+
C\|X(s)-\widetilde{X}_h(s)\|^2\Big(1+\|\widetilde{X}_h(s)\|_V^4+\|X(s)\|_V^4\Big)
\\
&
\quad
\leq
\tfrac12|\widetilde{e}_h(s)|_1^2
+
\tfrac98|\widetilde{e}_h(s)|_{-1,h}^2
+
C\|X(s)-\widetilde{X}_h(s)\|^2\Big(1+\|A_h\widetilde{X}_h(s)\|^4+|X(s)|_2^4\Big).
\end{split}
\end{align}
Integrating over $[0, t]$ and then using Gronwall's inequality one can arrive at
\begin{align}\label{eq:error-widetild-e-semi-H1-integrand}
|\widetilde{e}_h(t)|_{-1,h}^2
+
\int_0^t|\widetilde{e}_h(s)|_1^2\,\dd s
\leq
C
\int_0^t\|X(s)-\widetilde{X}_h(s)\|^2(1+\|A_h\widetilde{X}_h(s)\|^4+|X(s)|_2^4)\,\dd s
.
\end{align}
In view of \eqref{bound-AF}, \eqref{eq:bound-Ph-H2}, \eqref{lem:eq-smooth-property-Eh(t)}, 
\eqref{lem:eq-integrand-smooth-property-Eh(t)-II}
 and the Burkholder-Davis-Gundy inequality, we acquire that
\begin{align}\label{eq:bound-solution-semi-auxiliary}
\begin{split}
\|A_h\widetilde{X}_h(t)\|_{L^p(\Omega;\dot{H})}
\leq&
\|A_hE_h(t)P_hX_0\|_{L^p(\Omega;\dot{H})}
+
\int_0^t\|E_h(t-s)A^2_hP_h F(X(s))\|_{L^p(\Omega;\dot{H})}\,\dd s
\\
&+
C\Big(\int_0^t\|A_hE_h(t-s)P_hQ^{\frac12}\|^2_{\mathcal{L}_2}\,\dd s\Big)^{\frac12}
\\
\leq&
\|A_hP_hX_0\|_{L^p(\Omega;\dot{H})}
+
C\int_0^t(t-s)^{-\frac12}\|A_hP_h F ( X(s) )\|_{L^p(\Omega;\dot{H})}\,\dd s
+
C\|Q^{\frac12}\|_{\mathcal{L}_2}
\\
\leq&
C
\Big(
\|X_0\|_{L^p(\Omega;\dot{H}^2)}
+
\int_0^t(t-s)^{-\frac12}\,\dd s
\sup_{s\in[0,T]}\|PF(X(s))\|_{L^p(\Omega;\dot{H}^2)}
+
\|Q^{\frac12}\|_{\mathcal{L}_2}
\Big)
<\infty.
\end{split}
\end{align}
By employing \eqref{eq:error-auxilary-semi},
\eqref{eq:bound-solution-semi-auxiliary} and \eqref{them:spatial-regularity-mild-solution}, 
we derive from \eqref{eq:error-widetild-e-semi-H1-integrand} that
\begin{align}\label{eq:error-e-H1-integrand}
\begin{split}
\Big\|\int_0^t|&\widetilde{e}_h(s)|_1^2\,\dd s\Big\|_{L^p(\Omega;\mathbb{R})}
\leq
C\Big(
\int_0^t\Big\|\|X(s)-\widetilde{X}_h(s)\|^2(1+\|A_h\widetilde{X}_h(s)\|^4+|X(s)|_2^4
)\Big\|_{L^p(\Omega;\mathbb{R})}\dd s
\Big)
\\
&\leq
C\Big(\int_0^t\|X(s)-\widetilde{X}_h(s)\|_{L^{4p}(\Omega;\dot{H})}^4\,\dd s\Big)^{\frac12}
\Big(\int_0^t(1+\|A_h\widetilde{X}_h(s)\|_{L^{8p}(\Omega;\dot{H})}^8+\|X(s)\|_{L^{8p}(\Omega;\dot{H}^2)}^8)\,\dd s\Big)^{\frac12}
\\
&\leq
Ch^{2\kappa}|\ln h|^{2},
\qquad
\kappa=\min\{\gamma,r\}.
\end{split}
\end{align}
Equipped with this, we are ready  to bound
$\|\widetilde{e}_h(t)\|_{L^p(\Omega;\dot{H})}$, which can be decomposed by the following two terms:
\begin{align}\label{eq:error-e-H-decomposement}
\begin{split}
\|\widetilde{e}_h(t)\|_{L^p(\Omega;\dot{H})}
=&
\Big\|\int_0^tE_h(t-s)A_hP_hP\big(F(X(s))-F(X_h(s))\big)\,\dd s\Big\|_{L^p(\Omega;\dot{H})}
\\
\leq
&
\int_0^t\|E_h(t-s)A_hP_hP\big(F(X(s))-F(\widetilde{X}_h(s))\big)\|_{L^p(\Omega;\dot{H})}\,\dd s
\\
&+
\Big\|\int_0^tE_h(t-s)A_hP_hP\big(F(\widetilde{X}_h(s))-F(X_h(s))\big)\,\dd s\Big\|_{L^p(\Omega;\dot{H})}
\\
=:
&
J_1+J_2.
\end{split}
\end{align}
Following the same arguments as the proof of \eqref{eq:error-e-H1-integrand}, one can show
\begin{align}\label{eq:bound-L1-term}
\begin{split}
&J_1
\leq
C\int_0^t(t-s)^{-\frac12}\|F(\widetilde{X}_h(s))-F(X(s))\|_{L^p(\Omega;H)}\,\dd s
\\
&\leq C
\int_0^t(t-s)^{-\frac12}\|X(s)-\widetilde{X}_h(s)\|_{L^{2p}(\Omega;\dot{H})}
\big(1+\|A_h\widetilde{X}_h(s)\|_{L^{4p}(\Omega;\dot{H})}^2+\|X(s)\|_{L^{4p}(\Omega;\dot{H}^2)}^2\big)\,\dd s
\\
&\leq
Ch^\kappa|\ln h|,
\qquad
\kappa=\min\{\gamma,r\}.
\end{split}
\end{align}
Before handling the term $J_2$, we first adapt similar arguments used in the proof of \eqref{eq:bound-F-H1-semi}   and also use \eqref{eq:V-norm-control-by-Ah-norm} to get
\begin{align}\label{eq:bound-F-F-H1-semi}
\begin{split}
\|A_h^{\frac12}P_hP\big(F(X_h(s))-F(\widetilde{X}_h(s))\big)\|
\leq
&
\|A^{\frac12}P\big(F(X_h(s))-F(\widetilde{X}_h(s))\big)\|
\\
\leq
&
C|\widetilde{e}_h(s)|_1(1+\|A_hX_h(s)\|^2+\|A_h\widetilde{X}_h(s)\|^2).
\end{split}
\end{align}
This combined with \eqref{eq:error-e-H1-integrand}, \eqref{eq:bound-solution-semi-auxiliary}
and \eqref{lem:eq-bound-solution-semi} yields
\begin{align}\label{eq:bound-L2-term}
\begin{split}
J_2
&\leq
\Big\|
\int_0^t(t-s)^{-\frac14}
\big\|
A_h^{\frac12}P_h(F(X_h(s) ) - F ( \widetilde{X}_h(s)))
\big\|
\,\dd s
\Big\|_{L^p(\Omega;\mathbb{R})}
\\
&
\leq
C\Big\|\int_0^t(t-s)^{-\frac14}|\widetilde{e}_h(s)|_1
(1+\|A_h\widetilde{X}_h(s)\|^2+\|A_hX_h(s)\|^2)\,\dd s\Big\|_{L^p(\Omega;\mathbb{R})}
\\
&\leq
C\Big\|\Big(\int_0^t|\widetilde{e}_h(s)|_1^2\,\dd s\Big)^{\frac 12}
\Big(\int_0^t(t-s)^{-\frac12}\big(1+\|A_h\widetilde{X}_h(s)\|^2+\|A_hX_h(s)\|^2 \big)^2\,\dd s\Big)^{\frac12}\Big\|_{L^p(\Omega;\mathbb{R})}
\\
&\leq
C\Big\|\int_0^t|\widetilde{e}_h(s)|_1^2\,\dd s\Big\|_{L^p(\Omega;\mathbb{R})}^{\frac12}
\Big\|\int_0^t(t-s)^{-\frac12}\big(1+\|A_h\widetilde{X}_h(s)\|^2
+
\|A_hX_h(s)\|^2\big)^2\,\dd s\Big\|_{L^p(\Omega;\mathbb{R})}^{\frac1{2}}
\\
&\leq
Ch^\kappa|\ln h|,
\qquad
\kappa=\min\{\gamma,r\}.
\end{split}
\end{align}
Therefore, gathering the estimates of $J_1$ and $J_2$ together gives
\begin{align*}
\|\widetilde{X}_h(t)-X_h(t)\|_{L^p(\Omega;\dot{H})}
\leq
Ch^\kappa |\ln h|,
\end{align*}
which combined with \eqref{eq:error-auxilary-semi} validates \eqref{them:eq-error-estimates-semi-problem}.

We are now in the position to verify \eqref{them:eq-chemical potential-error-estimates-semi-problem}.
Similarly as before, we bound two terms $\|Y(t)-\widetilde{Y}_h(t)\|_{L^p(\Omega;\dot{H})}$
 and $\|\widetilde{Y}_h(t)-Y_h(t)\|_{L^p(\Omega;\dot{H})}$, where $\widetilde{Y}_h(t) := A_h\widetilde{X}_h(t)+P_hPF(X(t))$.
By \eqref{them:eq-mild-solution-stochastic-equation} and \eqref{eq:mild-soluiton-auxiliary-semi},
the error $\|Y(t)-\widetilde{Y}_h(t)\|_{L^p(\Omega;\dot{H})}$ can be decomposed as follows:
\begin{align}\label{eq:semi-auxiliary-error-potential}
\begin{split}
\|Y(t)-\widetilde{Y}_h(t)\|_{L^p(\Omega;\dot{H})}
\leq&
\underbrace{
\|(I-P_h)PF(X(t))\|_{L^p(\Omega;\dot{H})}
+
\|\big(AE(t)-A_hE_h(t)P_h\big)X_0\|_{L^p(\Omega;\dot{H})}
}_{ L_1 }
\\
&+
\underbrace{
\Big\|\int_0^t\big(A^2E(t-s)-A^2_hE_h(t-s)P_h)PF(X(s)\big)\,\dd s \Big\|_{L^p(\Omega;\dot{H})}
}_{L_2}
\\
&+
\underbrace{
\Big\|\int_0^t(AE(t-s)-A_hE_h(t-s)P_h)\,\dd W(s)\Big\|_{L^p(\Omega;\dot{H})}
}_{ L_3 }.
\end{split}
\end{align}
Using \eqref{eq:error-interpolation}, \eqref{bound-AF} and
\eqref{lem:eq-deterministic-error-displiace-semi} with $\alpha=2$ gives
\begin{align}\label{eq:L1-estimation}
L_1
\leq
C h^2 \sup_{s\in[0,T]}\|PF(X(s))\|_{L^p(\Omega;\dot{H}^2)}
+
Ch^2 t^{-1} \|X_0\|_{L^p(\Omega;\dot{H})}
\leq
Ch^2(1+t^{-1}).
\end{align}
To deal with the term $L_2$, we use \eqref{eq:relation-A-Ah-Rh-Ph} and the definition of the operator $\Phi_h(t)$ in \eqref{eq:definition-Psih-Phih(t)} to get
\begin{align}\label{eq:L2-decompose}
\begin{split}
L_2
\leq
&
\Big\|\int_0^t\Phi_h(t-s)APF(X(s))\,\dd s\Big\|_{L^p(\Omega;\dot{H})}
\\
&+
\Big\|\int_0^tA^2_hE_h(t-s)P_h(R_h-I)PF(X(s))\,\dd s\Big\|_{L^p(\Omega;\dot{H})}
\\
=:
&
L_{21}+L_{22}.
\end{split}
\end{align}
Owing to \eqref{bound-AF}, \eqref{lem:error-F(X)-F(Y)} with $\beta=2$, \eqref{lem:eq-chemical-potial-integrand-II} with $\varrho=2$ and
\eqref{lem:eq-deterministic-error-displiace-semi} with $\alpha=2$, we infer
\begin{align}\label{eq:L21-estimation}
\begin{split}
L_{21}
\leq
&
\Big\|\int_0^t\Phi_h(t-s)APF(X(t))\,\dd s\Big\|_{L^p(\Omega;\dot{H})}
+
\int_0^t\|\Phi_h(t-s)AP\big(F(X(s))-F(X(t))\big)\|_{L^p(\Omega;\dot{H})}\,\dd s
\\
\leq
&
Ch^2\|PF(X(t)\big)\|_{L^p(\Omega;\dot{H}^2)}
+
Ch^2\int_0^t(t-s)^{-1}\|P\big(F(X(s))-F(X(t))\big)\|_{L^p(\Omega;\dot{H}^2)}\,\dd s
\\
\leq
&
Ch^2\sup_{s\in[0,T]}\|PF(X(s))\|_{L^p(\Omega;\dot{H}^2)}
+
Ch^2\int_0^t(t-s)^{-1}(t-s)^{\frac14}\,\dd s
\\
\leq
&
Ch^2.
\end{split}
\end{align}
Likewise, we use  \eqref{bound-AF}, \eqref{lem:eq-integrand-smooth-property-Eh(t)}, \eqref{lem:error-F(X)-F(Y)} with $\beta=2$
and  \eqref{lem:eq-smooth-property-Eh(t)} with $\mu=2$ to derive
\begin{align}\label{eq:L22-estimation}
\begin{split}
L_{22}
\leq
&
\int_0^t\big\|A^2_hE_h(t-s)P_h(R_h-I)P\big(F(X(s))-F(X(t))\big)\big\|_{L^p(\Omega;\dot{H})}\,\dd s
\\
&
+
\Big\|\int_0^tA^2_hE_h(t-s)P_h(R_h-I)PF(X(t))\,\dd s\Big\|_{L^p(\Omega;\dot{H})}
\\
\leq
&
C h^2 \int_0^t(t-s)^{-1}  \|P\big(F(X(s))-F(X(t))\big)\|_{L^p(\Omega;\dot{H}^2)}\,\dd s
\\
&+
C\|( R_h - I ) P F( X( t ) ) \|_{L^p(\Omega; \dot{ H } )}
\\
\leq
&
Ch^2\int_0^t(t-s)^{-1}(t-s)^{\frac14}\,\dd s
+
Ch^2\sup_{s\in[0,T]}\| P F(X(s))\|_{L^p(\Omega;\dot{H}^2)}
\\
\leq
&
Ch^2,
\end{split}
\end{align}
which together with \eqref{eq:L21-estimation} implies
\begin{align}\label{eq:estimate-L22}
\begin{split}
L_2
\leq
Ch^2.
\end{split}
\end{align}
With the aid of  the Burkholder-Davis-Gundy inequality, one can use  \eqref{lem:eq-chemical-potial-integrand} with $\mu=\gamma-2$ to deduce
\begin{align}\label{eq:estimation-L3}
\begin{split}
L_3
\leq
C \Big(\int_0^t\big\|\Phi_h(t-s)Q^{\frac12}\big\|_{\mathcal{L}_2}^2\,\dd s\Big)^{\frac12}
\leq
C h^{\gamma-2} |\ln h| \|A^{\frac{\gamma-2}2}Q^{\frac12}\|_{\mathcal{L}_2}
\leq
Ch^{\gamma-2}|\ln h|.
\end{split}
\end{align}
Gathering \eqref{eq:L1-estimation}, \eqref{eq:estimate-L22} and \eqref{eq:estimation-L3} together implies
\begin{align}\label{eq:error-Y-widetild(Y)}
\|Y(t)-\widetilde{Y}_h(t)\|_{L^p(\Omega;\dot{H})}
\leq
Ch^{\gamma-2}|\ln h|(1+t^{-1}),
\qquad
\forall t \in (0, T].
\end{align}
To bound the error $\|\widetilde{Y}_h ( t_n ) - Y_h(t_n)\|_{L^p(\Omega;\dot{H})}$,
we first apply \eqref{eq:local-condition}, \eqref{eq:V-norm-control-by-Ah-norm},
\eqref{eq:embedding-equatlity-I}, \eqref{them:spatial-regularity-mild-solution}, \eqref{lem:eq-bound-solution-semi} 
and \eqref{them:eq-error-estimates-semi-problem} to achieve
\begin{align}\label{eq:error-FX-FXh}
\begin{split}
\|P\big(F(X(t))-F(X_h(t))\big)\|_{L^p(\Omega; \dot{H} ) }
\leq&
(1 + \|X( t )\|^2_{L^{4p}(\Omega;V)}
+
\|X_h(t)\|^2_{L^{4p}(\Omega;V)}
)
\|X(t)-X_h(t)\|_{L^{2p}(\Omega;\dot{H})}
\\
\leq
&
(1+\|X(t)\|^2_{L^{4p}(\Omega;\dot{H}^2)}
+
\|A_hX_h(t)\|^2_{L^{4p}(\Omega;\dot{H})}
)
\|X(t)-X_h(t)\|_{L^{2p}(\Omega;\dot{H})}
\\
\leq
&
Ch^\kappa |\ln h|.
\end{split}
\end{align}
Combining this  with  the inverse inequality \eqref{eq:inverse-ineq} enables us to obtain
\begin{align*}
\begin{split}
\|\widetilde{Y}_h(t_n)-Y_h(t_n)\|_{L^p(\Omega;\dot{H})}
\leq
&
\int_0^t\|A^2_hE_h(t-s)P_hP\big(F(X(s))-F(X_h(s))\big)\|_{L^p(\Omega;\dot{H})}\,\dd s
\\
&+
\big\| P_h P \big ( F(X(t))-F(X_h(t)) \big) \big\|_{L^p(\Omega;H)}
\\
\leq
&
Ch^{-1}\int_0^t(t-s)^{-\frac34}\|F(X(s))-F(X_h(s))\|_{L^p(\Omega;H)}\,\dd s
+
C h^\kappa|\ln h|
\\
\leq
&
C h^{\kappa-1}|\ln h|,
 \end{split}
\end{align*}
which in conjunction with \eqref{eq:error-Y-widetild(Y)} implies the desired assertion
\eqref{them:eq-chemical potential-error-estimates-semi-problem}.
$\square$
\section{The finite element full-discretization and its moment bounds}
In this section, we proceed to look at a finite element full discretization of \eqref{eq:abstract-SCHE} and provide
moment bounds of the full discretization. Let $k=T/N$, $N\in \mathbb{N}$ be a uniform time step-size
and $t_n=kn$, $n=1,2,\cdots, N$. We discretize \eqref{eq:semi-discrete-problem} in time with a 
backward Euler scheme and the resulting fully discrete problem is to  find $\mathcal{F}_{t_n}$-adapted
$\dot{V}_h$-valued random variable  $X_h^n $ such that
\begin{align}\label{eq:eq-solution-full}
X_h^n-X_h^{n-1}
+
kA_h^2X^n_h
+
k A_hP_hF(X_h^n)
=P_h\Delta W_n,
\; X_h^0=P_hX_0,
\quad
n=1,2,\cdots, N,
\end{align}
where $\Delta W_n := W(t_n)-W(t_{n-1})$.
Noting that the above implicit scheme works on the finite dimensional space $\dot{V}_h$ and that the mapping
$ A_h^2 + k A_h P_h F( \cdot ) $ obeys a kind of monotonicity condition in the Hilbert space $( \dot{V}_h,  ( \cdot , \cdot)_{ -1 , h } )$,
one can check that the implicit scheme \eqref{eq:eq-solution-full} is well-posed in $\dot{V}_h$.
After introducing a
family of operators  $\{E^n_{k,h}\}_{n=1}^N$:
\begin{align*}
E^n_{k,h}v_h
:=
(I + k A_h^2 )^{-n} v_h
=
\sum_{j=0}^{ \mathcal{N}_h } ( 1 + k \lambda_{j,h}^2 )^{-n} (v_h,e_{j,h})e_{j,h},
\quad
\forall v_h\in V_h,
\end{align*}
the solution of \eqref{eq:eq-solution-full}, similarly to the semi-discrete case,  
can be expressed in the following form,
\begin{align*}
X_h^n=E_{k,h}^nP_hX_0-k\sum_{j=1}^nA_hE_{k,h}^{n-j+1}P_hPF(X_h^j)+\sum_{j=1}^nE_{k,h}^{n-j+1}P_h\Delta W_j,
\quad
n=0, 1,2,\cdots, N.
\end{align*}
The next theorem offers a priori moment bounds for the fully-discrete approximations.
\begin{theorem}\label{lem:bound-numerical-solution}
Let $ X_h^n $ be given by \eqref{eq:eq-solution-full}.  If Assumptions
\ref{assum:linear-operator-A}-\ref{assum:intial-value-data}, \ref{assumption:triangulations} are valid, then there
 exist $k_0 > 0$ such that for all $k\leq k_0$, $h>0$ and $p\geq 1$,
\begin{align}\label{lem:eq-bound-solution-full-stochatic}
\sup_{1\leq n\leq N}\|A_hX^n_h\|_{L^p(\Omega;\dot{H})}
+
\Big\|\sum_{j=1}^Nk|A_hX_h^j+P_hF(X_h^j)|_1^2\Big\|_{L^p(\Omega; \R ) }
<\infty.
\end{align}
\end{theorem}
To prove it, we first introduce some smooth properties of the operator $E_{k,h}^n$.
Denoting $r(z) := (1+z)^{-1}$, one can write $E^n_{k,h}=r(kA_h^2)^n$.
As shown in \cite[Theorem 7.1]{ee2006galerkin}, there exist two constants $C$  and $c$ such that
\begin{align}
|r(z)-e^{-z}|
\leq
Cz^2,
\quad
\forall
z\in[0,1],
\end{align}
and
\begin{align}\label{eq:bound-R(z)-small}
|r(z)|
\leq e^{-cz},\;\forall z\in [0,1].
\end{align}
These two inequalities suffice to  ensure  that, for $n=1,2,3,\cdots,$
\begin{align}\label{eq:ez-r(z)-error}
|r(z)^n-e^{-zn}|
\leq
\Big|
\big(r(z)-e^{-z}\big)\sum_{l=0}^{n-1} r(z)^{n-1-l} e^{-zl}
\Big|
\leq
 C nz^2e^{-c(n-1)z},
 \quad
\forall
 z\in[0,1].
\end{align}

Additionally, we need a temporal discrete analogue of Lemma \ref{lem:propertiy-semgroup-semi} as follows.
\begin{lemma}\label{lem:bound-rational-approximation-semi}
Under Assumptions \ref{assum:linear-operator-A}, \ref{assumption:triangulations}, the following estimates hold for $n= 1, 2,3,  \cdots, N$,
\begin{align}
\label{lem:bounded-time-full-deterministic-problem}
\|A_h^\mu E_{k,h}^nP_hv\|
&\leq
Ct_n^{-\frac\mu2}\|v\|,
\quad \forall
\mu\in[0,2],
v \in \dot{H},
\\
\label{lem:temporal-reguarlity-time-full-deterministic-problem}
\|A_h^{-\nu}(I-E_{k,h}^n)P_hv\|
& \leq
C t_n^{\frac\nu2}\|v\|,\quad\forall \nu\in[0,2], v\in \dot{H},
\\
\label{lem:eq-bound-sum-full-operator}
\Big\|k\sum_{j=1}^nA_h^2 E_{k,h}^jP_h v\Big\|
& \leq
C\|v\|, \quad\forall v\in \dot{H},
\\
\label{lem:eq-bound-sum-full-operator-II}
\Big(k\sum_{j=1}^n\|A_h E_{k,h}^jP_h v\|^2\Big)^{\frac12}
& \leq
C\|v\|,\quad \forall v\in \dot{H}.
\end{align}
\end{lemma}
{\it Proof of Lemma \ref{lem:bound-rational-approximation-semi}.}
%
The estimate \eqref{lem:bounded-time-full-deterministic-problem} can be found in \cite[(2.10)]{furihata2018strong}.
And the proof of \eqref{lem:temporal-reguarlity-time-full-deterministic-problem} is based on an interpolation argument.
Taking $\mu=0$ in \eqref{lem:bounded-time-full-deterministic-problem} gives 
\eqref{lem:temporal-reguarlity-time-full-deterministic-problem}  for the case $\nu=0$. To show 
\eqref{lem:temporal-reguarlity-time-full-deterministic-problem}  for the case $\nu=2$,
we expand $P_hv$ in terms of $\big\{e_{j,h}\big\}_{j=1}^{\mathcal{N}_h}$ to obtain
\begin{align}\label{eq:temporal-reg-full}
\begin{split}
\|A_h^{-2}(I-E_{k,h}^n)P_hv\|^2
& =
\Big\|\sum_{j=1}^{\mathcal{N}_h}\lambda_{j,h}^{-2}\big(1-r(k\lambda_{j,h}^2)^n\big)(v,e_{j,h})e_{j,h}\Big\|^2
\\
&=
\sum_{j=1}^{\mathcal{N}_h}\lambda_{j,h}^{-4}\big(1-r(k\lambda^2_{j,h})^{n}\big)^2(v, e_{j,h})^2,
\end{split}
\end{align}
where $\big\{\lambda_{j,h}\big\}_{j=1}^{\mathcal{N}_h}$ are the positive eigenvalues of $A_h$ with corresponding
orthonormal eigenvectors $\big\{e_{j,h}\big\}_{j=1}^{\mathcal{N}_h}\subset \dot{V}_h$.
Since
\begin{align}
|1-r(k\lambda^2_{j,h})^{n}|
=
\big|1- (1+k\lambda^2_{j,h})^{-n} \big|
\leq
t_n \lambda^2_{j,h},
\quad
j=1,2,\cdots, \mathcal{N}_h,
\end{align}
one can derive from \eqref{eq:temporal-reg-full}  that \eqref{lem:temporal-reguarlity-time-full-deterministic-problem} holds for the case $\nu=2$.
 The intermediate cases follow by an interpolation argument.
To show \eqref{lem:eq-bound-sum-full-operator}, we again use Parseval's identity to infer
\begin{align*}
\Big\|k\sum_{j=1}^nA_h^2 E_{k,h}^jP_h v\Big\|^2
=
\Big\|k\sum_{j=1}^n\sum_{i=1}^{\mathcal{N}_h}\lambda_{i,h}^2\, r(k\lambda_{i,h}^2)^j(v,e_{i,h})e_{i,h}\Big\|^2
=
\sum_{i=1}^{\mathcal{N}_h} \Big(k\sum_{j=1}^n\lambda_{i,h}^2\, r(k\lambda_{i,h}^2)^j\Big)^2(v,e_{i,h})^2.
\end{align*}
Therefore,  \eqref{lem:eq-bound-sum-full-operator} holds on the condition,
\begin{align}\label{eq:summation-rkAh}
k\sum_{j=1}^n\lambda_{i,h}^2\, r(k\lambda_{i,h}^2)^j
\leq
 C,
 \quad
 i=1,2,\cdots, \mathcal{N}_h.
\end{align}
Next, we validate \eqref{eq:summation-rkAh} by considering two possibilities: either $k \lambda^2_{i,h} \leq 1$ or $k \lambda^2_{i,h} > 1$.
For the first possibility that $k \lambda^2_{i,h} \leq 1$, we use \eqref{eq:bound-R(z)-small} to deduce
\begin{align*}
\begin{split}
k\sum_{j=1}^n\lambda_{i,h}^2\, r(k\lambda_{i,h}^2)^j
\leq
k\sum_{j=1}^n\lambda_{i,h}^2\, e^{-ct_j\lambda_{i,h}^2}
\leq
\int_0^{t_n}\lambda_{i,h}^2\, e^{-cs\lambda_{i,h}^2}\,\dd s
\leq
\frac1 c(1-e^{-ct_n\lambda_{i,h}^2})
\leq
\frac 1 c.
\end{split}
\end{align*}
For the other possibility that $k \lambda^2_{i,h} > 1$, we notice that $r(k\lambda_{i,h}^2) < \frac12$ and thus
\begin{align*}
\begin{split}
k\sum_{j=1}^n\lambda_{i,h}^2\, r(k\lambda_{i,h}^2)^j
\leq
\sum_{j=1}^n r(k\lambda_{i,h}^2)^{j-1}
<
\sum_{j=1}^n2^{-(j-1)}
\leq
2.
\end{split}
\end{align*}
As a consequence, we obtain \eqref{eq:summation-rkAh} and the assertion \eqref{lem:eq-bound-sum-full-operator} follows.
The proof of \eqref{lem:eq-bound-sum-full-operator-II}  is similar to that of \eqref{lem:eq-bound-sum-full-operator} and we omit it.
$\square$

Equipped with the above preparations, we are prepared to show Theorem \ref{lem:bound-numerical-solution}.

{\it Proof of Theorem \ref{lem:bound-numerical-solution}.}
Following almost the  same lines as in the proof of \cite[Theorem 4.3]{furihata2018strong}, one can arrive at
\begin{align*}
\begin{split}
\Big\|
&\sup_{1\leq j\leq N} J(X_h^j)
\Big\|^p_{L^p(\Omega;\mathbb{R})}
+
\Big\|\sum_{j=1}^Nk|A_hX_h^j+P_hF(X_h^j)|_1^2\Big\|_{L^p(\Omega;\mathbb{R})}^p
\\
&\leq
C
\Big(
1+ \|J(P_hX_0)\|^p_{L^p(\Omega;\mathbb{R})}+\|P_hX_0\|^{4p+1}_{L^{4p+1}(\Omega;\dot{H})}
+
\|Q^{\frac12}(A_hP_hX_0+P_hF(P_hX_0)\,)\|_{L^{2p}(\Omega;\dot{H})}^{2p}
\Big)
\\
& <
\infty.
\end{split}
\end{align*}
The boundedness holds because $\|J(P_hX_0)\|^p_{L^p(\Omega;\mathbb{R})}+\|P_hX_0\|^{4p+1}_{L^{4p+1}(\Omega;\dot{H})}<\infty$ and
\begin{align*}
\begin{split}
\|Q^{\frac12}&(A_hP_hX_0+P_hF(P_hX_0)\,)\|_{L^{2p}(\Omega;\dot{H})}^{2p}
\\
\leq
&
C\|Q^{\frac12}\|_{\mathcal{L}(\dot{H})}^{2p}(\|A_hP_hX_0\|_{L^{2p}(\Omega;\dot{H})}^{2p}+\|P_hX_0\|_{L^{2p}(\Omega;\dot{H})}^{2p}+\|P_hX_0\|_{L^{6p}(\Omega;\dot{H}^1)}^{6p})
\\
\leq&
C(1+  \|X_0\|_{L^{2p}(\Omega;\dot{H}^2)}^{2p}
+
\|X_0\|_{L^{6p}(\Omega;\dot{H}^1)}^{6p})
<\infty,
\end{split}
\end{align*}
due to the use of \eqref{eq:embedding-equatlity-I}, \eqref{eq:bound-Ph-L4}, \eqref{eq:bound-Ph-H1} and \eqref{eq:bound-Ph-H2}.
Now it remains to bound the term $\|A_hX^n_h\|_{L^p(\Omega;\dot{H})}$.
Similarly to \eqref{eq:bound-semidiscrete-mild-H1}, we recall \eqref{eq:relation-H1-J} and obtain
\begin{align*}
\sup_{1\leq j\leq N}\|X_h^j\|_{L^p(\Omega;\dot{H}^1)}
\leq
C\Big\|\sup_{1\leq j\leq N}J(X_h^j)^{\frac p2}\Big\|_{L^1(\Omega;\mathbb{R})}^{\frac1{p}}
<
\infty,
\end{align*}
which can be used to promise
\begin{align}\label{eq:bound-F-full}
\sup_{1\leq j\leq N}\|F(X_h^j)\|_{L^p(\Omega;H)}
\leq
\Big(
\sup_{1\leq j\leq N}\|X_h^j\|_{L^p(\Omega;\dot{H})}
+
\sup_{1\leq j\leq N}\|X_h^j\|_{L^{3p} ( \Omega;\dot{H}^1 ) }^{3}
\Big)
<\infty.
\end{align}
With this at hand, one can follow the same lines in the proof of \eqref{eq:H2-bound-mild-solution-smei} to get
 \begin{align*}
 \sup_{1\leq n\leq N}\|A_hX_h^n\|_{L^p(\Omega;\dot{H})}<\infty.
 \end{align*}
The proof of this theorem is complete. $\square$
\section{Strong convergence rates of the FEM full-discretization}

The goal of this section is to identify strong convergence rates of
the fully discrete finite element method \eqref{eq:eq-solution-full}.
Similarly to the semi-discrete case, error estimates for the deterministic error operators 
and the moment bounds of the fully discrete finite element solutions
together play a key role in the convergence analysis.
To begin with,  we define the corresponding error operators
\begin{align}
\Psi_{k,h}(t):
=
E(t)-E_{k,h}^nP_h\;\; \text{and} \;\; \Phi_{k,h}(t):=AE(t)-A_hE_{k,h}^nP_h,\;\;t\in[t_{n-1},t_{n}),\;n\in\{1,2,\cdots,N\}.
\end{align}
Since the constant eigenmodes are canceled, it is easy to see
\begin{equation}
\label{eq:Psi-h-Phi-h-full}
\Phi_h(t) v = \Phi_h(t) P v,
\quad
v \in H.
\end{equation}
The following lemma is a temporal version of Lemma \ref{lem:deterministic-error-semi}, which is crucial in the error analysis.
\begin{lemma}\label{lem:error-estimes-time-full-deterministic-problem}
Under Assumptions \ref{assum:linear-operator-A}, \ref{assumption:triangulations}, 
the following estimates for $\Psi_{k,h}(t)$ and $\Phi_{k,h}(t)$ hold
\begin{align}
\label{lem:error-stimates-nonsmooth}
\|\Psi_{k,h}(t) v\|
&\leq
C(h^\beta+k^{\frac \beta 4})  |v|_\beta,
\quad
\forall \beta\in[1,r], \: v \in \dot{H}^\beta,
\: t\in[0,T],
\\
\|\Phi_{k,h}(t) v\|
&\leq
C(h^\alpha+k^{\frac \alpha 4}) t^{-1} |v|_{\alpha-2},
\quad
\forall \alpha\in[1,r], \: v\in \dot{H}^{\alpha-2}, t \in (0, T],
\label{lem:eq-error-error-estimes-time-full-deterministic-problem}
\\
\label{lem:error-F-full-integrand}
\Big(\int_0^t\|\Psi_{k,h}(s)v\|^2\, \mathrm{d} s\Big)^{\frac12}
& \leq
C(h^{\nu}|\ln h|+ |\ln k| k^{\frac\nu4})|v|_{\nu-2},
\quad
\forall v\in \dot{H}^{\nu-2},
\: \nu\in[1,r], \: t \in [ 0, T],
\\
\label{lem:error-deterministc-potial-full-integrand}
\Big(\int_0^t \|\Phi_{k,h}(s)v\|^2\,\mathrm{d} s\Big)^{\frac12}
& \leq
C(h^{\mu}|\ln h|+ k^{\frac\mu4}|\ln k|)|v|_\mu,
\quad
\forall \mu\in[0,r],
v\in \dot{H}^{\mu},
t \in [ 0, T],
\\
\label{lem:error-deterministc-potial-full-integrand-II}
\Big\|\int_0^t \Phi_{k,h}(s)v \,\mathrm{d} s\Big\|
& \leq
C(h^\varrho+k^{\frac\varrho4})|v|_{\varrho-2},
\quad
\forall \varrho\in[1,r], \: v\in \dot{H}^{\varrho-2},
t \in [0, T],
\end{align}
where $r\in\{2,3,4\}$ for $d=1$ and $r=2$ for $d\in\{2,3\}$, as implied by Assumption \ref{assumption:triangulations}.
\end{lemma}

The proof of Lemma \ref{lem:error-estimes-time-full-deterministic-problem}  is postponed to Appendix. 
Equipped with this lemma we are well-prepared to do the error analysis.
The next theorem states the main result of this section, concerning strong convergence rates of the FEM full-discretization.
\begin{theorem}\label{them:error-estimates-full-stochastic}
Let $X(t)$ be the weak solution of \eqref{eq:abstract-SCHE} and let  $X_h^n$ be given by \eqref{eq:eq-solution-full}. 
Let Assumptions \ref{assum:linear-operator-A}-\ref{assum:intial-value-data}
be valid for some $\gamma\in[3,4]$ and let Assumption \ref{assumption:triangulations} be fulfilled with 
$r\in\{2,3,4\}$ for $d=1$ and $r=2$ for $d\in\{2,3\}$. 
Then there exist $k_0 > 0$ such that for all $k\leq k_0$, $h>0$, $p\geq 1$ and $ n \in \{ 0, 1,..., N\}$,
\begin{align}\label{them:error-estimates-full-problem}
\|X(t_n)-X_h^n\|_{L^p(\Omega;\dot{H})}
\leq
C(h^\kappa|\ln h|+k^{\frac\kappa4} |\ln k|)
\quad
\text{ with }
\;
\kappa : = \min\{\gamma,r\}.
\end{align}
Moreover, the discrepancy between the "chemical potential" $Y(t) :=  AX(t)+PF(X(t))$ and its approximation $Y_h^n := A_h X_h^n+P_hPF(X^n_h)$ is measured as follows, for any $k\leq k_0$, $h>0$, $p\geq 1$ and $n \in \{ 1,2,..., N\}$,
\begin{align}\label{them:eq-chemical-potential-error-full-estimates-problem}
\|Y(t_n)-Y_h^n\|_{L^p(\Omega;\dot{H})}
\leq
C(1+t_n^{-1}\big)(h^\iota |\ln h|+k^{\frac\iota4}|\ln k|)
\quad
\text{ with }
\;
\iota=\min\{\gamma-2,r-1\}.
\end{align}
\end{theorem}
%
{\it Proof of Theorem \ref{them:error-estimates-full-stochastic}.}
Similarly to the semi-discrete case, by introducing the auxiliary problem,
\begin{align}\label{eq:auxiliary-problem}
\widetilde{X}_h^n-\widetilde{X}_h^{n-1}
+
k A_h\big(A_h\widetilde{X}_h^n+P_hF(X(t_n))\big)
=
P_h\Delta W_n,
\quad 
\widetilde{X}_h^0 = P_h X_0,
\end{align}
whose solution can be recasted as
\begin{align}\label{eq:solution-form-auxiliary-problem}
\widetilde{X}_h^n
=
E_{k,h}^n P_h X_0
-
k\sum_{j=1}^n A_hE_{k,h}^{n-j+1}P_hF(X(t_j))
+
\sum_{j=1}^n E_{k,h}^{n-j+1} P_h \Delta W_j,
\end{align}
we decompose the considered error $\|X(t_n)-X_h^n\|_{L^p(\Omega;\dot{H})}$ into two parts:
\begin{align}
\|X(t_n)-X_h^n\|_{L^p(\Omega;\dot{H})}
\leq
\|X(t_n)-\widetilde{X}_h^n\|_{L^p(\Omega;\dot{H})}
+
\|\widetilde{X}_h^n-X_h^n\|_{L^p(\Omega;\dot{H})}.
\end{align}
At first we handle the estimate of the first error term. Subtracting \eqref{eq:solution-form-auxiliary-problem} from
\eqref{them:eq-mild-solution-stochastic-equation}, we split it into the following three parts:
\begin{align}\label{eq:full-auxiliary-error}
\begin{split}
\|X(t_n)-\widetilde{X}_h^n\|_{L^p(\Omega;\dot{H})}
\leq &
\|(E(t_n)-E_{k,h}^nP_h)X_0\|_{L^p(\Omega;\dot{H})}
\\
&+
\Big\|\int_0^{t_n}E(t_n-s)APF(X(s))\,\dd s
-
k\sum_{j=1}^nA_hE_{k,h}^{n-j+1}P_hPF((X(t_j))\Big\|_{L^p(\Omega;\dot{H})}
\\
&+
\Big\|\sum_{j=1}^n\int_{t_{j-1}}^{t_j}\big(E(t_n-s)
-
E_{k,h}^{n-j+1}P_h\big) \,\dd W(s)\Big\|_{L^p(\Omega;\dot{H})}
\\
=:& \,
\mathbb{I}+\mathbb{J}+\mathbb{K}.
\end{split}
\end{align}
The three terms $\mathbb{I}, \mathbb{J}, \mathbb{K}$ will be treated separately.
By using \eqref{lem:error-stimates-nonsmooth} with $\beta = \kappa$,
we estimate the first term $\mathbb{I}$ as follows,
\begin{align*}
\mathbb{I}
\leq
C(h^\kappa+k^{\frac\kappa4})\|X_0\|_{L^p(\Omega;\dot{H}^\kappa)},
\quad
\kappa = \min\{\gamma,r\}.
\end{align*}
To bound the term $\mathbb{J}$,  we need to decompose it further:
\begin{align*}
\begin{split}
\mathbb{J}
\leq
&
\bigg\|
\sum_{j=1}^n\int_{t_{j-1}}^{t_j}  E(t_n-s)AP\big(F(X(s))-F(X(t_j))\big)\,\dd s
\bigg\|_{L^p(\Omega;\dot{H})}
\\
&
+
\Big\|\int_0^{t_n}\Phi_{k,h}(t_n-s)PF(X(t_n))\,\dd s\Big\|_{L^p(\Omega;\dot{H})}
\\
&+
\sum_{j=1}^
n\int_{t_{j-1}}^{t_j}\|\Phi_{k,h}(t_n-s)P\big(F(X(t_j))-F(X(t_n))\big)\|_{L^p(\Omega;\dot{H})}\,\dd s
\\
=:
&
\mathbb{J}_1+\mathbb{J}_2+\mathbb{J}_3.
\end{split}
\end{align*}
Noticing that, for $s\in[t_{j-1},t_j)$,
\begin{align*}
X(t_j)
=
E(t_j-s)X(s)
-
\int_{s}^{t_j}E(t_j-\sigma)APF(X(\sigma))\,\dd \sigma
+
\int_{s}^{t_j}E(t_j-\sigma)\,\dd W(\sigma),
\end{align*}
and using Taylor's formula help us to split $\mathbb{J}_1$ into four additional terms:
\begin{align*}
\begin{split}
\mathbb{J}_1
\leq
&
\Big\|\sum_{j=1}^n\int_{t_{j-1}}^{t_j}E(t_n-s)APF'(X(s))(E(t_j-s)-I)X(s)\,\dd s\Big\|_{L^p(\Omega;\dot{H})}
\\
&+
\Big\|\sum_{j=1}^n\int_{t_{j-1}}^{t_j}E(t_n-s)APF'(X(s))\int_s^{t_j}E(t_j-\sigma)A P F(X(\sigma))\,\dd \sigma\,\dd s\Big\|_{L^p(\Omega;\dot{H})}
\\
&+
\Big\|\sum_{j=1}^n\int_{t_{j-1}}^{t_j}E(t_n-s)APF'(X(s))\int_s^{t_j}E(t_j-\sigma)\,\dd W(\sigma)\,\dd s \Big\|_{L^p(\Omega;\dot{H})}
\\
&+
\Big\|\sum_{j=1}^n\int_{t_{j-1}}^{t_j}E(t_n-s)APR_F(X(s),X(t_j))\,\dd s\Big\|_{L^p(\Omega;\dot{H})}
\\
= :
&
\,
\mathbb{J}_{11}+\mathbb{J}_{12}+\mathbb{J}_{13}+\mathbb{J}_{14}.
\end{split}
\end{align*}
Here the remainder term $R_F$ is given by,
\begin{align*}
R_F(X(s),X(t_j))
:=
\int_0^1
F''(X(s)+\lambda(X(t_j)-X(s)))
\big( X(t_j)-X(s), X(t_j)-X(s) \big) ( 1 - \lambda )\,\dd \lambda.
\end{align*}
In view of \eqref{I-spatial-temporal-S(t)} with $\nu=\frac{2+\delta_0}{2}$, \eqref{eq:one-seoncd-derivation-f}, 
\eqref{them:spatial-regularity-mild-solution}, \eqref{eq:embedding-equatlity-III} and the H\"{o}lder inequality,  
we derive, for any fixed $\delta_0\in(\frac32,2)$ and $\gamma \in[3,4]$,
\begin{align}\label{eq:estimate-L2111}
\begin{split}
\mathbb{J}_{11}
&\leq
C\sum_{j=1}^n\int_{t_{j-1}}^{t_j}(t_n-s)^{-\frac{2+\delta_0}4}
      \|A^{-\frac{\delta_0}2}PF'(X(s))(E(t_j-s)-I)X(s)\|_{L^{p}(\Omega;\dot{H})}\,\dd s
\\
&\leq
C\sum_{j=1}^n\int_{t_{j-1}}^{t_j}(t_n-s)^{-\frac{2+\delta_0}4}\|F'(X(s))
    (E(t_j-s)-I)X(s)\|_{L^{p}(\Omega;L_1)}\,\dd s
\\
&\leq
C\sum_{j=1}^n\int_{t_{j-1}}^{t_j}(t_n-s)^{-\frac{2+\delta_0}4}
      \|f'(X(s))\|_{L^{2p}(\Omega;H)}\|(E(t_j-s)-I)X(s)\|_{L^{2p}(\Omega;\dot{H})}\,\dd s
\\
&\leq
Ck^{\frac\gamma4}\int_0^{t_n}(t_n-s)^{-\frac{2+\delta_0}4}\,\dd s
\sup_{s\in[0,T]}\|f'(X(s))\|_{L^{2p}(\Omega;H)}\sup_{s\in[0,T]}\|X(s)\|_{L^{2p}(\Omega;\dot{H}^\gamma)}
\\
&\leq
 Ck^{\frac\gamma4}.
\end{split}
\end{align}
Following similar arguments as in \eqref{eq:estimate-L2111}, 
we use \eqref{bound-AF}  to show that, for any $ \delta_0 \in ( \tfrac32, 2 ), $
\begin{align}\label{eq:estimate-L2112}
\begin{split}
\mathbb{J}_{12}
&\leq
\sum_{j=1}^n\int_{t_{j-1}}^{t_j}\int_{s}^{t_j}(t_n-s)^{-\frac{2+\delta_0}4}
\big\|F'(X(s))E(t_{i+1}-\sigma)A P F(X(\sigma))\big\|_{L^{p}(\Omega;L_1)}\,\dd \sigma\,\dd s
\\
&\leq
C\sum_{j=1}^n\int_{t_{j-1}}^{t_j}\int_s^{t_j}(t_n-s)^{-\frac{2+\delta_0}4}
\|f'(X(s))\|_{L^{2p}(\Omega;H)}\|PF(X(\sigma))\|_{L^{2p}(\Omega;\dot{H}^2)}\,\dd \sigma\,\dd s
\\
&\leq
Ck \int_0^{t_n}(t_n-s)^{-\frac{2+\delta_0}4} \dd s \sup_{s\in[0,T]}\|f'(X(s))\|_{L^{2p}(\Omega;H)}
\sup_{s\in[0,T]}\|PF(X(s))\|_{L^{2p}(\Omega;\dot{H}^2)}
\\
&\leq
 Ck.
\end{split}
\end{align}
When estimating $\mathbb{J}_{13}$, we recall the stochastic Fubini theorem (see \cite[Theorem 4.18]{da2014stochastic})
and the Burkholder-Davis-Gundy-type inequality to obtain
\begin{align*}
\begin{split}
\mathbb{J}_{13}
=&
\Big\|\sum_{j=1}^n\int_{t_{j-1}}^{t_j}\int_{t_{j-1}}^{t_j}\chi_{[s,t_j)}(\sigma)
   E(t_n-s)APF'(X(s))E(t_j-\sigma)\,\dd W(\sigma)\,\dd s \Big\|_{L^{p}(\Omega;\dot{H})}
\\
=
&
\Big\|\sum_{j=1}^n\int_{t_{j-1}}^{t_j}\int_{t_{j-1}}^{t_j}\chi_{[s,t_j)}(\sigma)
  E(t_n-s)APF'(X(s))E(t_j-\sigma)\,\dd s \dd W(\sigma) \Big\|_{L^{p}(\Omega;\dot{H})}
\\
\leq
&
\Big(\sum_{j=1}^n\int_{t_{j-1}}^{t_j}\Big\|\int_{t_{j-1}}^{t_j}
\chi_{[s,t_j)}(\sigma)E(t_n-s)APF'(X(s))E(t_j-\sigma)Q^{\frac12}\,\dd s\Big\|_{L^p(\Omega;\mathcal{L}_2)}^2\,\dd \sigma\Big)^{\frac12},
\end{split}
\end{align*}
where $\chi_{[s,t_j)}(\cdot)$ stands for the indicator function defined by
$\chi_{[s,t_j)}( \sigma )
=
1
\text{ for } \sigma \in [s,t_j)
\text{ and }
\chi_{[s,t_j)}( \sigma ) = 0
\text{ for } \sigma \not \in [s,t_j).
$
Further, we employ the H\"{o}lder inequality, \eqref{lem:bound-A-deriveaton-f} and \eqref{III-spatial-temporal-S(t)} with $\varrho=\frac{4-\gamma}2$ to deduce,
\begin{align}\label{eq:j13-estimatation-full}
\begin{split}
\mathbb{J}_{13}
\leq
&
Ck^{\frac12}
\Big(\sum_{j=1}^n\int_{t_{j-1}}^{t_j}\int_{t_{j-1}}^{t_j}\sum_{l=1}^\infty\|
E(t_n-s)APF'(X(s))E(t_j-\sigma)Q^{\frac12} e_l\|_{L^p(\Omega;\dot{H})}^2\,\dd s\,\dd \sigma\Big)^{\frac12}
\\
\leq
&
Ck^{\frac12} \Big(\sum_{j=1}^n\int_{t_{j-1}}^{t_j}
\sum_{l=1}^\infty\int_{t_{j-1}}^{t_j}(t_n-s)^{-\frac12}\|A^{\frac12}PF'(X(s))E(t_j-\sigma)Q^{\frac12} e_l \|_{L^p(\Omega;\dot{H})}^2\,\dd s\dd \sigma\Big)^{\frac12}
\\
\leq
&
Ck^{\frac12}  \Big(\sum_{j=1}^n\int_{t_{j-1}}^{t_j}
\sum_{l=1}^\infty\int_{t_{j-1}}^{t_j}(t_n-s)^{-\frac12}
\Big(
1
+
\sup_{ r \in [ 0, T] }
\|
X ( r )
\|^2_{ L^{2p}(\Omega; \dot{H}^2 ) }
\Big)
\|A^{\frac12}E(t_j-\sigma)Q^{\frac12} e_l \|_{L_6}^2\,\dd s\dd \sigma\Big)^{\frac12}
\\
\leq
&C k^{\frac12}  \Big(\sum_{j=1}^n\int_{t_{j-1}}^{t_j}(t_n-s)^{-\frac12}
\,\dd s
\int_{t_{j-1}}^{t_j}\|AE(t_j-\sigma)Q^{\frac12}\|_{\mathcal{L}_2}^2\,\dd \sigma\Big)^{\frac12}
\\
\leq
&
Ck^{\frac\gamma4} \Big(\int_0^{t_n}(t_n-s)^{-\frac12}\,\dd s\Big)^{\frac12}
\|A^{\frac{\gamma-2}2}Q^{\frac12}\|_{\mathcal{L}_2}
\\
\leq
&
Ck^{\frac\gamma4},
\end{split}
\end{align}
where $\gamma\in[3,4]$ comes from the assumption \eqref{eq:ass-AQ-condition} and $\{ e_i\}_{i \in \N}$ is 
any orthogonal basis basis of $\dot{H}$.
To bound the term $\mathbb{J}_{14}$,  we use \eqref{them:temporal-regularity-mild-stoch}, \eqref{eq:embedding-equatlity-I}
and \eqref{eq:one-seoncd-derivation-f} to infer
\begin{align*}
\mathbb{J}_{14}
&\leq
C\sum_{j=1}^n\int_{t_{j-1}}^{t_j}(t_n-s)^{-\frac{2+\delta_0}4}
     \|R_F(X(s),X(t_j))\|_{L^p(\Omega;L_1)}\,\dd s
\nonumber\\
&
\leq
C\sum_{j=1}^n\int_{t_{j-1}}^{t_j}(t_n-s)^{-\frac{2+\delta_0}4}
\big\|\;\|X(t_j)-X(s)\|
\,\|f''((1-\lambda)X(s)+\lambda X(t_j))\|_{L_4}\,\|X(t_j)-X(s)\|_{L_4}\;\big\|_{L^p(\Omega;\mathbb{R})}
\,\dd s
\nonumber\\
&
\leq
C\sum_{j=1}^n\int_{t_{j-1}}^{t_j}(t_n-s)^{-\frac{2+\delta_0}4}\|X(t_j)-X(s)\|_{L^{4p}(\Omega;\dot{H})}\|X(t_j)-X(s)\|_{L^{4p}(\Omega;\dot{H}^1)}\,\dd s
\sup_{s\in[0,T]}\|f''(X(s))\|_{L^{2p}(\Omega;L_4)}
\nonumber\\
&\leq
Ck\sup_{s\in[0,T]}\|f''(X(s))\|_{L^{2p}(\Omega;L_4)}
         \int_0^{t_n}(t_n-s)^{-\frac{2+\delta_0}4}\,\dd s
\nonumber
\\
&
\leq Ck,
\end{align*}
where for the first step we followed similar arguments as used in \eqref{eq:estimate-L2111}.
This together with \eqref{eq:estimate-L2111}, \eqref{eq:estimate-L2112} and \eqref{eq:j13-estimatation-full}
leads to, for $\gamma\in[3,4]$,
\begin{align*}
\mathbb{J}_1
\leq
Ck^{\frac\gamma4}.
\end{align*}
Concerning  the term $\mathbb{J}_2$, we apply  \eqref{bound-AF} and \eqref{lem:error-deterministc-potial-full-integrand-II} with $\varrho=\kappa=\min\{\gamma,r\}$ to get
\begin{align*}
\begin{split}
\mathbb{J}_2
&\leq
C(h^\kappa + k^{\frac\kappa4} ) \|PF(X(t_n))\|_{L^p(\Omega;\dot{H}^{\kappa-2})}
\leq
C(h^\kappa + k^{\frac\kappa4} ) \|PF(X(t_n))\|_{L^p(\Omega;\dot{H}^2)}
\leq
C(h^\kappa + k^{\frac\kappa4} ).
\end{split}
\end{align*}
With regard to $\mathbb{J}_3$, after employing \eqref{lem:eq-error-error-estimes-time-full-deterministic-problem} with $\alpha=\kappa$
 and \eqref{lem:error-F(X)-F(Y)} with $\beta=2$ one can arrive at
\begin{align*}
\begin{split}
\mathbb{J}_3
&\leq
C(h^\kappa+ k^{\frac\kappa4} )\sum_{j=1}^
n\int_{t_{j-1}}^{t_j}(t_n-s)^{-1}\|A^{\frac{\kappa-2}2}P\big(F(X(t_j))-F(X(t_n))\big)\|_{L^p(\Omega;\dot{H})}\,\dd s
\\
&\leq
C(h^\kappa+ k^{\frac\kappa4} )\sum_{j=1}^
{n-1}\int_{t_{j-1}}^{t_j}(t_n-s)^{-1} t_{n-j}^{\frac14}\,\dd s
\\
&\leq
C(h^\kappa+ k^{\frac\kappa4} ),
\quad
\kappa = \min\{\gamma,r\},
\end{split}
\end{align*}
where we also used the facts $\kappa-2=\min\{r-2,\gamma-2\}\leq 2$ and $t_{n-j}^\frac14 \leq (t_n - s)^{\frac14}$ for $s \leq t_j$.
Gathering the above three estimates together results in
\begin{align*}
\mathbb{J}
\leq
C(h^\kappa + k^{\frac \kappa4}  ),
\quad
\kappa=\min\{\gamma,r\}.
\end{align*}
For the last term $\mathbb{K}$,
we  utilize \eqref{lem:error-F-full-integrand} with $\nu=\kappa$ and the Burkholder-Davis-Gundy inequality to obtain
\begin{align*}
\begin{split}
\mathbb{K}
\leq
\Big(\sum_{j=1}^n\int_{t_{j-1}}^{t_j}\|\Psi_{k,h}(t_n-s)Q^{\frac12}\|^2_{\mathcal{L}_2}\dd s\Big)^{\frac12}
\leq
C(h^\kappa|\ln h|+k^{\frac\kappa4} |\ln k|) \|A^{\frac{\kappa-2}2}Q^{\frac12}\|_{\mathcal{L}_2},
\,
\kappa=\min\{\gamma,r\}.
\end{split}
\end{align*}
Now, putting the above estimates together results in
\begin{align}\label{eq:x-x^n-estiamte}
\|X(t_n)-\widetilde{X}_h^n\|_{L^p(\Omega;\dot{H})}
\leq
C(h^\kappa|\ln h|+k^{\frac\kappa4}|\ln k|),
\quad
\kappa=\min\{\gamma,r\}.
\end{align}
Next we turn our attention to the error $\widetilde{e}_h^n := X_h^n-\widetilde{X}_h^n$, obeying
\begin{align}\label{eq:error-equaton-discrete}
\widetilde{e}_h^n-\widetilde{e}_h^{n-1}
+
kA_h^2\widetilde{e}_h^n=-kA_hP_hF(X_h^n)+kA_hP_hF(X(t_n)),\quad \widetilde{e}_h^0=0.
\end{align}
Equivalently, this can be reformulated as
\begin{align}\label{eq:solution-error-equation-full-discrete}
\widetilde{e}_h^n=k\sum_{j=1}^nA_hE_{k,h}^{n-j+1}P_h(F(X(t_j))-F(X_h^j)).
\end{align}
Before proceeding further, we need to bound the term $\|A_h\widetilde{X}_h^n\|_{L^p(\Omega;\dot{H})}$.
Owing to \eqref{bound-AF}, \eqref{lem:bounded-time-full-deterministic-problem}, \eqref{lem:eq-bound-sum-full-operator-II},
\eqref{eq:bound-Ph-H2}
 and the Burkholder-Davis-Gundy inequality, one can derive, for any $n\in\{1,2,\cdots,N\}$,
\begin{align}\label{eq:boundness-solution-form-full-auxiliary-problem}
\begin{split}
\|A_h\widetilde{X}_h^n\|_{L^p(\Omega;\dot{H})}
\leq&
\|A_hE_{k,h}^nP_hX_0\|_{L^p(\Omega;\dot{H})}
+
k\sum_{j=1}^n\|A^2_hE_{k,h}^{n-j+1}P_hPF(X(t_j))\|_{L^p(\Omega;\dot{H})}
\\
&
+
\Big\|\sum_{j=1}^nA_hE_{k,h}^{n-j+1}P_h\Delta W^j\Big\|_{L^p(\Omega;\dot{H})}
\\
\leq
&
C\|A_hP_hX_0\|_{L^p(\Omega;\dot{H})}
+
Ck\sum_{j=1}^nt_{n-j+1}^{-\frac12}\|A_hP_hPF(X(t_j))\|_{L^p(\Omega;\dot{H})}
\\
&
+
C\Big(k\sum_{j=1}^n\|A_hE_{k,h}^{n-j+1}P_hQ^{\frac12}\|_{\mathcal{L}_2}^2\Big)^{\frac12}
\\
\leq
&
C
\Big(
1+\|X_0\|_{L^p(\Omega;\dot{H}^2)}
+
k\sum_{j=1}^{n}t_{n-j+1}^{-\frac12}
\sup_{s\in[0,T]}\|PF(X(s))\|_{L^p(\Omega;\dot{H}^2)}
+
\|Q^{\frac12}\|_{\mathcal{L}_2}
\Big)
<\infty.
\end{split}
\end{align}
Multiplying both sides of \eqref{eq:error-equaton-discrete} by $A_h^{-1}\widetilde{e}_h^n$ yields
\begin{align*}
(\widetilde{e}_h^n-\widetilde{e}_h^{n-1},A_h^{-1}\widetilde{e}_h^n)
+
k(\nabla \widetilde{e}_h^n,\nabla \widetilde{e}_h^n)
=
k(-F(X_h^n)+F(\widetilde{X}_h^n),\widetilde{e}_h^n)
+
k(-F(\widetilde{X}_h^n)+F(X(t_n)),\widetilde{e}_h^n).
\end{align*}
Noting that $\widetilde{e}_h^0=0$ and $\frac12(|\widetilde{e}_h^n|_{-1,h}^2-|\widetilde{e}_h^{n-1}|_{-1,h}^2)
\leq
(\widetilde{e}_h^n-\widetilde{e}_h^{n-1},A_h^{-1}\widetilde{e}_h^n)$,
one can follow a similar way as in \eqref{eq:error-widetild-e-semi} to arrive at
\begin{align*}
\begin{split}
\tfrac12(|\widetilde{e}_h^n|_{-1,h}^2-|\widetilde{e}_h^{n-1}|_{-1,h}^2)
+
k| \widetilde{e}_h^n|_1^2
\leq
\tfrac{k}2|\widetilde{e}_h^n|_1^2
+
\tfrac{9k}8|\widetilde{e}_h^n|_{-1,h}^2
+
Ck\|\widetilde{X}_h^n-X(t_n)\|^2
(
1 + |\widetilde{X}_h^n |_{2,h}^4
+
|X(t_n)|_2^4
)
.
\end{split}
\end{align*}
Summation on $n$ and applying the Gronwall inequality give
\begin{align*}
|\widetilde{e}_h^n|_{-1,h}^2
+
k\sum_{j=1}^n|\widetilde{e}_h^j|_1^2
\leq
Ck\sum_{j=1}^n \|\widetilde{X}_h^j-X(t_j)\|^2(1+|\widetilde{X}^j_h|_{2,h}^4+|X(t_j)|_2^4),
\end{align*}
which together with \eqref{them:spatial-regularity-mild-solution},
\eqref{eq:boundness-solution-form-full-auxiliary-problem} and \eqref{eq:x-x^n-estiamte} leads to
\begin{align}\label{eq:sum-x-x^n-in-H1}
\begin{split}
\Big\|k\sum_{j=1}^n|\widetilde{e}_h^j|_1^2\Big\|_{L^p(\Omega;\mathbb{R})}
&\leq
C
k\sum_{j=1}^n \Big\|\|\widetilde{X}_h^j-X(t_j)\|^2(1+|\widetilde{X}^j_h|_{2,h}^4+|X(t_j)|_2^4)\Big\|_{L^p(\Omega;\mathbb{R})}
\\
&\leq
Ck\sum_{j=1}^n \|\widetilde{X}^j_h-X(t_j)\|_{L^{4p}(\Omega; \dot{H})}^2
\big(1+\|A_h\widetilde{X}_h^j\|_{L^{8p}(\Omega;\dot{H})}^4+\|X(t_j)\|_{L^{8p}(\Omega;\dot{H}^2)}^4\big)
\\
&\leq
C(h^\kappa|\ln h|+k^{\frac\kappa4}|\ln k|)^2,
\quad
\kappa=\min\{\gamma,r\}.
\end{split}
\end{align}
Similarly to the semi-discrete case, we use \eqref{eq:solution-error-equation-full-discrete} to
split the error $\|\widetilde{e}_h^n\|_{L^p(\Omega; \dot{H})}$ as follows:
\begin{align}\label{eq:widetilde{e}-decompose-full}
\begin{split}
\|\widetilde{e}_h^n\|_{L^p(\Omega;\dot{H})}
\leq&
k\sum_{j=1}^n\|A_hE_{k,h}^{n-j+1} P_h P (F(X(t_j))-F(\widetilde{X}_h^j))\|_{L^p(\Omega;\dot{H})}
\\
&+
k\|\sum_{j=1}^nA_hE_{k,h}^{n-j+1} P_h P(F(\widetilde{X}_h^j)-F(X_h^j))\|_{L^p(\Omega;\dot{H})}
\\
=:
&
\mathbb{A}+\mathbb{B}.
\end{split}
\end{align}
Similarly to \eqref{eq:bound-L1-term}, we employ \eqref{eq:x-x^n-estiamte},
\eqref{lem:bounded-time-full-deterministic-problem}, \eqref{eq:boundness-solution-form-full-auxiliary-problem}
 and  \eqref{them:spatial-regularity-mild-solution}  to get
\begin{align}\label{eq:A-estimate-full}
\begin{split}
\mathbb{A}
&\leq
Ck\sum_{j=1}^nt_{n-j+1}^{-\frac12}\|F(X(t_j))-F(\widetilde{X}_h^j)\|_{L^p(\Omega; \dot{H})}
\\
&
\leq Ck\sum_{j=1}^nt_{n-j+1}^{-\frac12}\|X(t_j)-\widetilde{X}_h^j\|_{L^{2p}(\Omega;\dot{H})}(1+\|X(t_j)\|^2_{L^{4p}(\Omega;\dot{H}^2)}
+
\|A_h\widetilde{X}_h^j\|^2_{L^{4p}(\Omega;\dot{H})})
\\
&
\leq
C(h^\kappa|\ln h|+k^{\frac\kappa4}|\ln k |) k\sum_{j=1}^nt_{n-j+1}^{-\frac12} (1+\sup_{s\in[0,T]}\|X(s)\|^2_{L^{4p}(\Omega;\dot{H}^2)}
+
\sup_{1\leq j\leq N}\|A_h\widetilde{X}_h^j\|^2_{L^{4p}(\Omega;\dot{H})})
\\
&
\leq
C(h^\kappa| \ln h| + k^{\frac\kappa4}| \ln k|),
\quad
\kappa=\min\{\gamma,r\}.
\end{split}
\end{align}
For the term $\mathbb{B}$, similar techniques used in \eqref{eq:bound-F-F-H1-semi} help us to show
\begin{align}\label{eq:bound-I-full}
\|A_h^{\frac12}P_hP(F(\widetilde{X}_h^j)-F(X_h^j))\|
\leq
C|\widetilde{X}_h^j-X_h^j|_1(1+\|A_h\widetilde{X}_h^j\|^2+\|A_hX_h^j\|^2).
\end{align}
Combining this with \eqref{eq:sum-x-x^n-in-H1}, \eqref{eq:boundness-solution-form-full-auxiliary-problem}
and \eqref{lem:eq-bound-solution-full-stochatic} enables us to derive
\begin{align*}
\begin{split}
\mathbb{B}
\leq
&
\Big\|k\sum_{j=1}^nt_{n-j+1}^{-\frac14}\|A_h^{\frac12}P_hP(F(\widetilde{X}_h^j)-F(X_h^j))\|\Big\|_{L^p(\Omega;\mathbb{R})}
\\
\leq
&
C\Big\|k\sum_{j=1}^nt_{n-j+1}^{-\frac14}|\widetilde{e}_h^j|_1
\,
(1+\|A_h\widetilde{X}_h^j\|^2+\|A_hX_h^j\|^2)\Big\|_{L^p(\Omega;\mathbb{R})}
\\
\leq
&
C\Big\|\Big(k\sum_{j=1}^n |\widetilde{e}_h^j|^2_1\Big)^{\frac 12}
\Big(k\sum_{j=1}^nt_{n-j+1}^{-\frac12}(1+\|A_h\widetilde{X}_h^j\|^4+\|A_hX_h^j\|^4)\Big)^{\frac 12}\Big\|_{L^p(\Omega;\mathbb{R})}
\\
\leq
&
C\Big\|k\sum_{j=1}^n |\widetilde{e}^j_h|_1^2\Big\|_{L^p(\Omega;\mathbb{R})}^{\frac12}
\Big\|k\sum_{j=1}^nt_{n-j+1}^{-\frac12}(1+\|A_h\widetilde{X}_h^j\|^4+\|A_hX_h^j\|^4)\Big\|_{L^p(\Omega;\mathbb{R})}^{\frac12}
\\
\leq
&
C(h^\kappa|\ln h|+k^{\frac\kappa4}|\ln k|),
\quad
\kappa=\min\{\gamma,r\},
\end{split}
\end{align*}
which together with \eqref{eq:A-estimate-full}, \eqref{eq:widetilde{e}-decompose-full} and \eqref{eq:x-x^n-estiamte}
gives the desired assertion \eqref{them:error-estimates-full-problem}.

In the sequel, we focus on the error $\|Y(t_n)-Y_h^n\|_{L^p(\Omega; \dot{H})}$.
Similarly to the semi-discrete case, we first consider the error $\|Y(t_n)-\widetilde{Y}_h^n\|_{L^p(\Omega;\dot{H})}$,
 where $ \widetilde{Y}_h^n = A_h\widetilde{X}_h^n + P_hPF(X(t_n))$. By \eqref{eq:solution-form-auxiliary-problem}
 and \eqref{them:eq-mild-solution-stochastic-equation},
\begin{align}\label{eq:full-auxiliary-chemical-error}
\begin{split}
\|Y(t_n)-\widetilde{Y}_h^n\|_{L^p(\Omega;\dot{H})}
\leq
&
\underbrace{
\|(AE(t_n)-A_hE_{k,h}^nP_h)X_0\|_{L^p(\Omega;\dot{H})}
+
\|(I-P_h)PF(X(t_n))\|_{L^p(\Omega;\dot{H})}
}_{ \mathbb{L}_1 }
\\
&
+
\underbrace{
\Big\|\sum_{j=1}^n\int_{t_{j-1}}^{t_j}A^2E(t_n-s)PF(X(s))-A_h^2E_{k,h}^{n-j+1}P_hPF(X(t_j))\,\dd s\Big\|_{L^p(\Omega;\dot{H})}
}_{\mathbb{L}_2}
\\
&
+
\underbrace{
\Big(\sum_{j=1}^n\int_{t_{j-1}}^{t_j}\|(AE(t_n-s)-A_hE^{n-j+1}_{k,h}P_h)Q^{\frac12}\|_{\mathcal{L}_2}^2\,\dd s\Big)^{\frac12}
}_{ \mathbb{L}_3 }.
\end{split}
\end{align}
In the same spirit as in \eqref{eq:L1-estimation} but employing
\eqref{lem:eq-error-error-estimes-time-full-deterministic-problem} with $\alpha=2$ instead we obtain
\begin{align}\label{eq:estimation-L1}
\begin{split}
\mathbb{L}_1
&\leq
Ch^2\sup_{s\in[0,T]}\|PF(X(s))\|_{L^p(\Omega;\dot{H}^2)}
+
C(h^2+k^{\frac12})t_n^{-1} \|X_0\|_{L^p(\Omega;\dot{H})}
\\
&\leq
C(h^2+k^{\frac12})(1+t_n^{-1}).
\end{split}
\end{align}
In order to properly handle $\mathbb{L}_2$, we need its further decomposition  as follows,
\begin{align*}
\begin{split}
\mathbb{L}_2
\leq
&
\Big\|\sum_{j=1}^n\int_{t_{j-1}}^{t_j}A^2E(t_n-s)P\big(F(X(s))-F(X(t_j))\big)\,\dd s\Big\|_{L^p(\Omega;\dot{H})}
\\
&+
\sum_{j=1}^n\int_{t_{j-1}}^{t_j}\|(A^2E(t_n-s)-A^2_hE_{k,h}^{n-j+1}P_h)P\big(F(X(t_j)-F(X(t_n))\,\big)\|_{L^p(\Omega;\dot{H})}\,\dd s
\\
&+
\Big\|\sum_{j=1}^n\int_{t_{j-1}}^{t_j}\big(A^2E(t_n-s)-A^2_hE_{k,h}^{n-j+1}P_h\big)PF(X(t_n))\,\dd s\Big\|_{L^p(\Omega;\dot{H})}
\\
=: &
\mathbb{L}_{21}+\mathbb{L}_{22}+\mathbb{L}_{23}.
\end{split}
\end{align*}
Thanks to
\eqref{lem:error-F(X)-F(Y)} with $\beta=1$  and \eqref{I-spatial-temporal-S(t)} with $\mu=\frac32$,
one can show
\begin{align}\label{eq:l1-estimation-full}
\begin{split}
\mathbb{L}_{21}
& \leq
C\sum_{j=1}^n\int_{t_{j-1}}^{t_j}(t_n-s)^{-\frac34}\|A^{\frac12}P\big(F(X(s))-F(X(t_j))\big)\|_{L^p(\Omega;\dot{H})}\,\dd s
\\
& \leq
Ck^{\frac12} \int_0^{t_n}(t_n-s)^{-\frac34}\,\dd s
\leq
Ck^{\frac12}.
\end{split}
\end{align}
Similarly to \eqref{eq:L2-decompose}, using \eqref{lem:error-F(X)-F(Y)} with $\beta=2$, \eqref{lem:eq-error-error-estimes-time-full-deterministic-problem} with $\alpha=2$, \eqref{lem:bounded-time-full-deterministic-problem} with $\mu=2$  and  \eqref{eq:error-interpolation} implies
\begin{align*}
\begin{split}
\mathbb{L}_{22}
\leq
&
\sum_{j=1}^{n-1}\int_{t_{j-1}}^{t_j}\|\Phi_{k,h}(t_n-s)AP\big(F(X(t_j))-F(X(t_n))\big)\|_{L^p(\Omega;\dot{H})}\,\dd s
\\
&
+
\sum_{j=1}^{n-1}\int_{t_{j-1}}^{t_j}\|A^2_hE_{k,h}^{n-j+1}P_h(I-R_h)P\big(F(X(t_j))-F(X(t_n))\big)\|_{L^p(\Omega;\dot{H})}\,\dd s
\\
\leq
&
C(h^2+k^{\frac12})\sum_{j=1}^{n-1}\int_{t_{j-1}}^{t_j}\big((t_n-s)^{-1}+t_{n-j+1}^{-1}\big)\|AP\big(F(X(t_j))-F(X(t_n))\big)\|_{L^p(\Omega;\dot{H})}\,\dd s
\\
\leq
&
C(h^2+k^{\frac12} ) \sum_{j=1}^{n-1} k(t_{n-j}^{-1+\frac14}+t_{n-j+1}^{-1+\frac14})
\\
\leq
&
C(h^2+k^{\frac12} ).
\end{split}
\end{align*}
Similarly as before, we utilize \eqref{eq:relation-A-Ah-Rh-Ph}, \eqref{lem:eq-bound-sum-full-operator}, \eqref{bound-AF},
  \eqref{lem:error-deterministc-potial-full-integrand-II} with $\varrho=2$
    and \eqref{eq:error-interpolation} to bound $\mathbb{L}_{23}$ as follows,
\begin{align*}
\mathbb{L}_{23}
\leq
&
\Big\|\int_0^{t_n}\Phi_{k,h}(t_n-s)APF(X(t_n))\,\dd s \Big\|_{L^p(\Omega;\dot{H})}
+
\Big\|\sum_{j=1}^nk A^2_hE_{k,h}^{n-j+1}P_h(R_h-I)PF(X(t_n)) \Big\|_{L^p(\Omega;\dot{H})}
\nonumber\\
\leq
&
C(h^2 + k^{\frac12} ) \|PF(X(t_n))\|_{L^p(\Omega;\dot{H}^2)}
+
C\|(R_h-I)PF(X(t_n))\|_{L^p(\Omega;\dot{H})}
\nonumber\\
\leq
&
C(h^2 + k^{\frac12} ) \sup_{s\in[0,T]}\|PF(X(s))\|_{L^p(\Omega;\dot{H}^2)}
.
\end{align*}
Putting the above three estimates together ensures
\begin{align}\label{eq:L_2-estimate-full}
\mathbb{L}_2
\leq
C(h^2 + k^{\frac12} ).
\end{align}
At the moment we start to estimate the term $\mathbb{L}_3$. In the light of
\eqref{lem:error-deterministc-potial-full-integrand} with $\mu = \gamma-2 $, we derive
\begin{align*}
\begin{split}
\mathbb{L}_3
\leq
\Big(\sum_{j=1}^n\int_{t_{j-1}}^{t_j}\|\Phi_{k,h}(t_n-s)Q^{\frac12}\|_{\mathcal{L}_2}^2\,\dd s\Big)^{\frac12}
\leq
C(h^{\gamma-2}|\ln h|+k^{\frac{\gamma-2}4}|\ln k|)\|A^{\frac{\gamma-2}2}Q^{\frac12}\|_{\mathcal{L}_2},
\end{split}
\end{align*}
which together with \eqref{eq:L_2-estimate-full} and \eqref{eq:estimation-L1} shows
\begin{align}\label{eq:error-Y-wide(Y)-full}
\|Y(t_n)-\widetilde{Y}_h^n\|_{L^p(\Omega;\dot{H})}
\leq
C(h^{\gamma-2} |\ln h| + k^{\frac{\gamma-2}4} |\ln k|  ) ( 1 + t_n^{-1}).
\end{align}
Now it remains to bound $\|\widetilde{Y}_h^n-Y_h^n\|_{L^p(\Omega;\dot{H})}$. 
Using the same arguments as in \eqref{eq:error-FX-FXh} promises
\begin{align*}
\begin{split}
&\|P\big(F(X(t_j))-F(X_h^j)\big)\|_{L^p(\Omega;\dot{H})}
\\
& \quad \leq
C (1+\sup_{s\in[0,T]}\|X(s)\|^2_{L^{4p}(\Omega;\dot{H}^2)}
+
\sup_{1\leq j\leq N} \|A_hX_h^j\|^2_{L^{4p}(\Omega;H)})
\|X(t_j)-X_h^j\|_{L^{2p}(\Omega;\dot{H})}
\\
& \quad \leq
C(h^\kappa|\ln h|+k^{\frac\kappa4} |\ln k|),
\quad
\kappa=\min\{\gamma,r\}.
\end{split}
\end{align*}
Combining this with \eqref{lem:bounded-time-full-deterministic-problem} with $ \mu = \tfrac12 $, the inverse inequality \eqref{eq:inverse-ineq}
and the fact $t^{-1}_{n-j+1}\leq Ck^{-1}$ helps us to arrive at
\begin{align}\label{eq:eror-widet(Y)-Y-full}
\begin{split}
\|\widetilde{Y}_h^n-Y_h^n\|_{L^p(\Omega;\dot{H})}
\leq
&
\sum_{j=1}^nk\|E_{k,h}^{n-j+1}A^2_hP_hP(F(X(t_j))-F(X_h^j))\|_{L^p(\Omega;\dot{H})}
\\
&
+
\|P\big(F(X(t_n))-F(X_h^n)\big)\|_{L^p(\Omega;\dot{H})}
\\
\leq
&
C\min\{h^{-1},k^{-\frac14}\} \sum_{j=1}^nkt_{n-j+1}^{-\frac34}\|P\big(F(X(t_j))-F(X_h^j)\big)\|_{L^p(\Omega;\dot{H})}
\\
&+
C(h^\kappa|\ln h|+ k^{\frac\kappa4}|\ln k|)
\\
\leq
&
C(h^{\kappa-1}|\ln h|+ k^{\frac{\kappa-1}4}|\ln k|),
\quad
\kappa=\min\{\gamma,r\}.
\end{split}
\end{align}
Putting \eqref{eq:error-Y-wide(Y)-full}  and \eqref{eq:eror-widet(Y)-Y-full} together
finally gives \eqref{them:eq-chemical-potential-error-full-estimates-problem}, as required.
$\square$
\section{Numerical experiments}
\label{sect:numerical-section}
Some numerical tests are presented in this section to illustrate the previous
findings.  We consider the following Cahn-Hilliard-Cook equation in one dimension,
\begin{align}\label{eq:numericla-expeirent-eq}
\left\{\begin{array}{ll}
\frac{\partial u}{\partial t}
=\frac{\partial^2w}{\partial x^2}+\dot{W},&t\in(0,T],\;x\in(0,1),\\
w=- \frac{\partial^2u}{\partial x^2}+ u-u^3,&x\in(0,1),\\
u(0,x)=\cos(\pi x),&x\in(0,1),\\
\frac{\partial u}{\partial x}|_{x=0}=\frac{\partial u}{\partial x}|_{x=1}=0,&t\in(0,T],\\
\frac{\partial w}{\partial x}|_{x=0}=\frac{\partial w}{\partial x}|_{x=1}=0,&t\in(0,T],
\end{array}
\right.
\end{align}
where $\{W(t)\}_{t\in[0,T]}$ is a standard $Q$-Wiener process in $\dot{H}$, with a covariance
operator $Q = A^{-1.5005}$. Here $-A$ is the Laplacian with homogeneous Neumann boundary conditions
and the condition \eqref{eq:ass-AQ-condition} is fulfilled with $\gamma = 3$.
We aim to perform mean-square approximations of the exact solution $(u, w)$ to \eqref{eq:numericla-expeirent-eq}
at the endpoint  $T = 1$. To do so we take a piecewise linear finite element method for the spatial discretization
and a backward Euler method for the temporal discretization.  
The expectations are approximated by computing averages over $M = 100$ samples.  
The exact solutions $(u, w)$, not available at hand, are computed by numerical ones with small step-sizes $h_{exact}$ and $k_{exact}$.

\begin{figure}[!ht]
\centering
      \includegraphics[width=3in,height=3in] {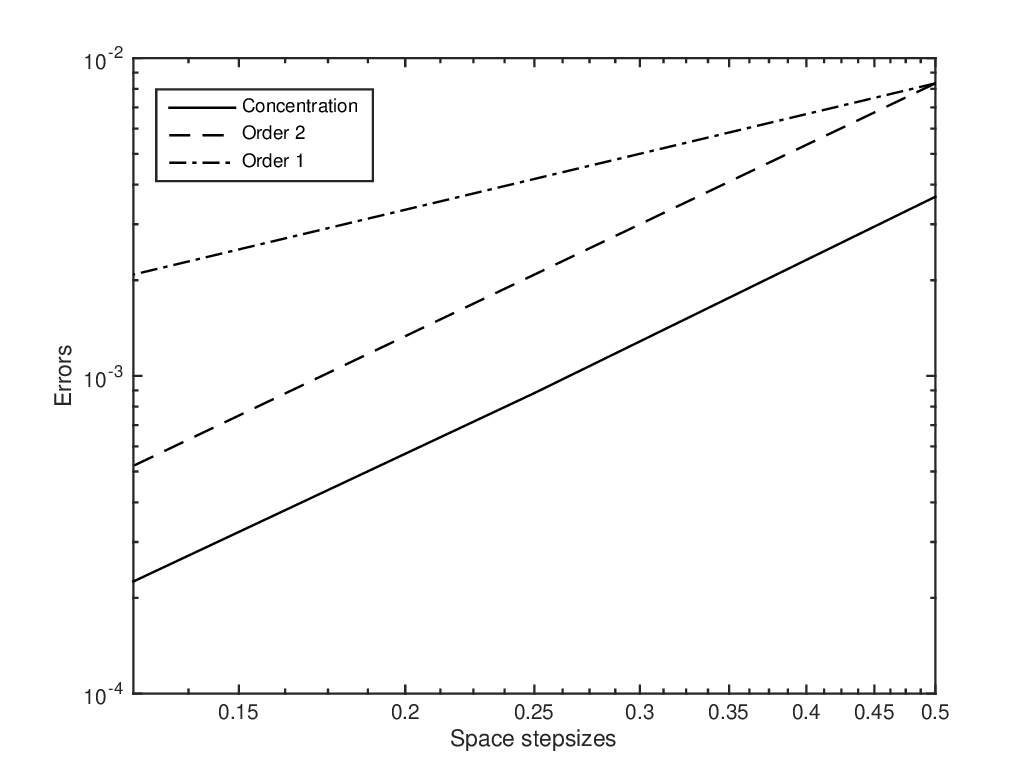}
       \includegraphics[width=3in,height=3in] {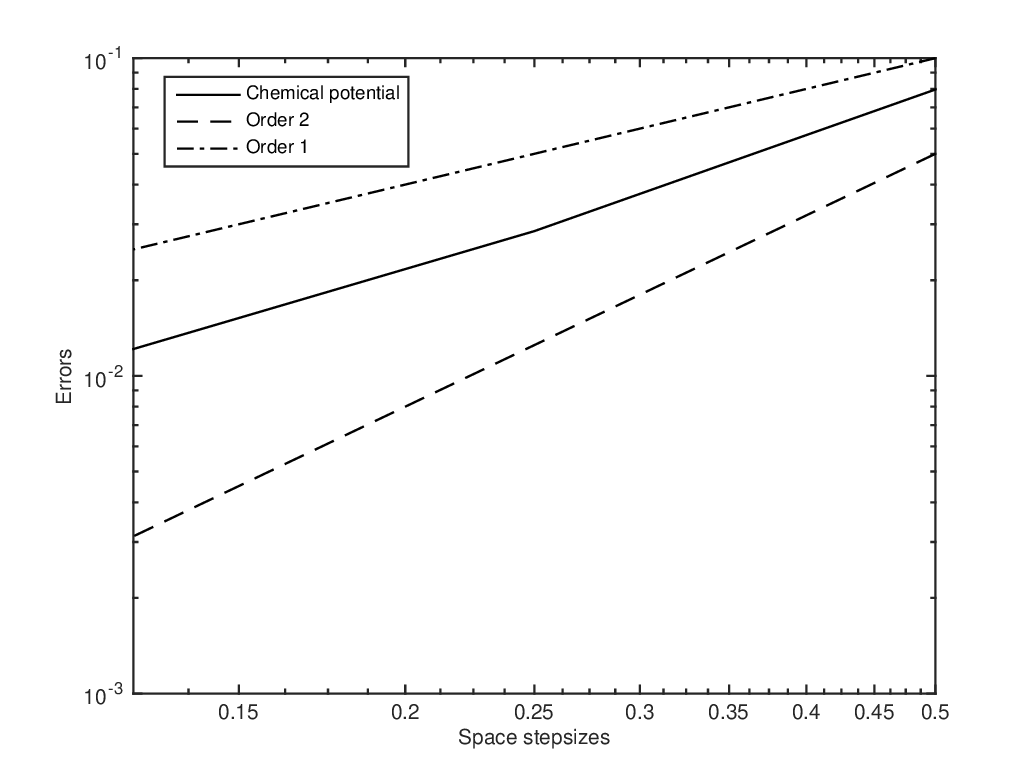}
 \caption{Mean-square convergence rates for the spatial discretizations}
 \label{picutre:space-numerical-experiment}
\end{figure}

\begin{figure}[!ht]
\centering
      \includegraphics[width=3in,height=3in] {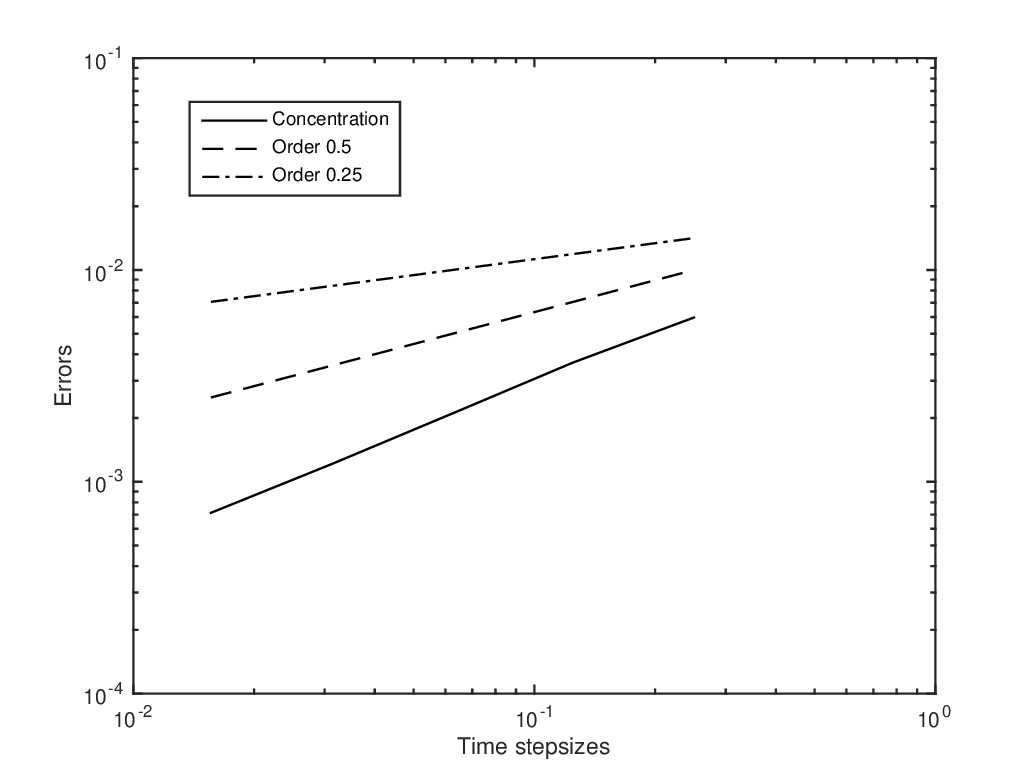}
       \includegraphics[width=3in,height=3in] {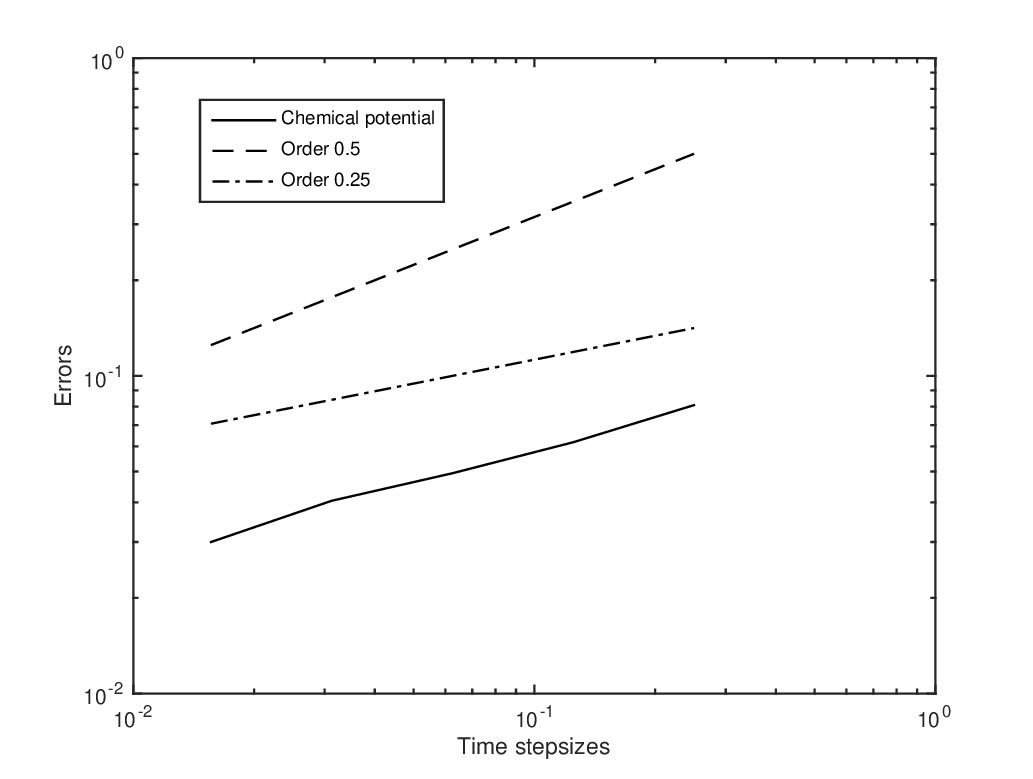}
 \caption{Mean-square convergence rates for the temporal discretizations }\label{picure:temporal-numerical-experiment}
\end{figure}

At first, we test  the convergence rates for the spatial discretizations.  The ``true" solutions $(u,w)$ are computed 
using  $h_{exact}=2^{-5}$ and $k_{exact}=2^{-20}$. 
We carry out numerical simulations with three different spatial step-sizes $h=2^{-i}, i=1,2,3$ and 
present the resulting mean-square errors (solid lines) in Figure \ref{picutre:space-numerical-experiment}. 
As expected, convergence rates of order $2$ and order $1$ are detected 
for the concentration and the chemical potential, respectively. 
This is consistent  with the previous theoretical findings in Theorem \ref{them:error-estimates-semi-problem} 
when $\gamma=3$ and $r=2$. 
%
Next we turn to the temporal discretizaton and fix $h=2^{-7}$ for the FEM semi-discretization,  
whose ``true" solutions $(u_h,w_h)$ are computed 
using  $k_{exact}=2^{-14}$. Similarly,  we present in Figure \ref{picure:temporal-numerical-experiment} 
errors due to the temporal discretizatons using five time stepsizes $k=2^{-i},i=6,7,8,9,10$.  
Comparing the slopes of the error (solid) lines with those of the reference (dashed) lines,  one can see,
temporal approximation errors decrease at a slope close to order $\frac12$
for the concentration and order $\frac14$ for the chemical potential,  which agree with the theoretical results.

%

\section{Appendix: Proof of Lemma \ref{lem:error-estimes-time-full-deterministic-problem}}
%
The estimates  \eqref{lem:error-stimates-nonsmooth} and \eqref{lem:error-F-full-integrand} can be proved
 by a simple modification of the proof of \cite[Theorem 2.2]{larsson2011finite}.
In order to validate \eqref{lem:eq-error-error-estimes-time-full-deterministic-problem}, one can first use  
\eqref{I-spatial-temporal-S(t)} with $\mu=2$, \eqref{II-spatial-temporal-S(t)} with $\nu=\frac\alpha2$ 
and \eqref{lem:eq-deterministic-error-displiace-semi} to get, for $ t \in [ t_{n-1},t_n )$,
 \begin{align}\label{eq:phiv}
 \begin{split}
 \|\Phi_{k,h}(t) v\|
 &\leq
 \|A(E(t)-E(t_n))v\|
 +
 \|(AE(t_n)-A_hE_h(t_n)P_h) v\|
 +
 \|A_h(E_h(t_n)-E_{k,h}^n)P_hv\|
 \\
& \leq
 \|A^2E( t )A^{-\frac\alpha2}(I-E(t_n - t )) A^{\frac{\alpha-2}2} v\|
 +
 Ct_n^{-1}h^\alpha |v|_{\alpha-2}
 +
  \|A_h(E_h(t_n)-E_{k,h}^n)P_hv\|
  \\
 &
 \leq
 Ct^{-1}(h^\alpha+k^{\frac\alpha4})|v|_{\alpha-2}
 +
  \|A_h(E_h(t_n)-E_{k,h}^n)P_hv\|,
  \end{split}
 \end{align}
where
 \begin{align}\label{eq:Ah(Eh-Ekh)}
  \|A_h(E_h(t_n)-E_{k,h}^n)P_hv\|
  \leq
  Ct_n^{-1}\|A_h^{-1}P_hv\|,
 \end{align}
 due to the use of \eqref{lem:eq-smooth-property-Eh(t)} with $\mu=2$ and \eqref{lem:bounded-time-full-deterministic-problem} with $\mu=2$.
On the other hand, \cite[Theorem 4.4]{larsson1991finite} shows
\begin{align*}
  \|A_h(E_h(t_n)-E_{k,h}^n)P_hv\|
  \leq
  Ckt_n^{-1}\|A_hP_hv\|.
 \end{align*}
An interpolation between these two estimates shows, for $\beta\in[0,4]$ and $t\in[t_{n-1},t_n)$,
\begin{align*}
\|A_h(E_h(t_n)-E_{k,h}^n)P_hv\|
\leq
Ct_n^{-1} k^{\frac\beta4}\|A_h^{\frac {\beta-2}2}P_hv\|
\leq
Ct^{-1} k^{\frac\beta4}\|A_h^{\frac {\beta-2}2}P_hv\|,
\end{align*}
which, after assigning $\beta=\alpha\in[1,r]$ and considering \eqref{eq:phiv} and \eqref{eq:Ah-A-bound}, implies \eqref{lem:eq-error-error-estimes-time-full-deterministic-problem}.
Repeating the same arguments in the proof of  \eqref{lem:eq-chemical-potial-integrand},
we can show \eqref{lem:error-deterministc-potial-full-integrand}.
Next we prove \eqref{lem:error-deterministc-potial-full-integrand-II}.  Note first that, for $t\in[t_n,t_{n+1})$, $n\geq 0$,
\begin{align}\label{eq:estimate-Phi(k,h)(s)-integrand-decompose}
\begin{split}
\Big\|\int_0^t \Phi_{k,h}(s)v \,\mathrm{d} s\Big\|
\leq
&
\Big\|\int_{t_n}^t\Phi_{k,h}(s)v\,\mathrm{d} s\Big\|
+
\Big\|\int_0^{t_n} \Phi_{k,h}(s)v \,\mathrm{d} s\Big\|.
\end{split}
\end{align}
By virtue of \eqref{eq:relation-A-Ah},  \eqref{II-spatial-temporal-S(t)} with $\nu=\frac\varrho2$ and \eqref{lem:bounded-time-full-deterministic-problem} with $\mu=\frac{4-\varrho}2$, we acquire
\begin{align}\label{eq:estimate-Phi(k,h)(s)-integrand-I}
\begin{split}
\Big\| \int_{t_n}^t \Phi_{k,h}(s)v\,\mathrm{d} s\Big\|
&=
\big\|A^{-1}\big(E(t)-E(t_n)\big)v\big\|
+
(t-t_n)\|A_hE_{k,h}^{n+1}P_hv\|
\\
&\leq
C
\|A^{-\frac\varrho2}\big(I-E(t-t_n)\big) A^{\frac{\varrho-2}2}v\|
+
(t-t_n)
\|A_h^{\frac{4-\varrho}2}E_{k,h}^{n+1} A_h^{\frac{\varrho-2}2}P_hv\|
\\
&\leq
C\big(k^{\frac\varrho4}
+
(t-t_n)t_{n+1}^{-\frac{4-\varrho}4}\big)|v|_{\varrho-2}
\\
&
\leq
Ck^{\frac\varrho4} |v|_{\varrho-2}.
\end{split}
\end{align}
Further, owing to \eqref{lem:eq-chemical-potial-integrand-II} we obtain
\begin{align}\label{eq:estimate-Phi(k,h)(s)-integrand-II}
\begin{split}
\Big\| \int_0^{t_n}\Phi_{k,h}(s)v\,\mathrm{d} s\Big\|
\leq &
\Big\|\int_0^{t_n}\Phi_{h}(s)v\,\mathrm{d} s\Big\|
+
\Big\|\sum_{j=1}^n\int_{t_{j-1}}^{t_j}\big( E_h(s)-E_h(t_j)\big)A_hP_hv\,\mathrm{d} s\Big\|
\\
&+
\Big\|\sum_{j=1}^n\int_{t_{j-1}}^{t_j}\big( E_h(t_j)-E_{k,h}^j\big)A_hP_hv\,\mathrm{d} s\Big\|
\\
\leq
& Ch^\varrho|v|_{\varrho-2}
+
\Big\|\sum_{j=1}^n\int_{t_{j-1}}^{t_j}\big( E_h(s)-E_h(t_j)\big)A_hP_hv\,\mathrm{d} s\Big\|
\\
& +
\Big\|k\sum_{j=1}^n\big( E_h(t_j)-E_{k,h}^j\big)A_hP_hv\,\Big\|,
\end{split}
\end{align}
where using Parseval's identity and the fact $\lambda^{-\frac\varrho2}_{i,h}(1-e^{-(t_j-s)\lambda_{i,h}^2})\leq C k^{\frac\varrho4}$, $0\leq\varrho\leq 4$  gives
\begin{align}\label{eq:summand-first-full}
\begin{split}
\Big\|\sum_{j=1}^n\int_{t_{j-1}}^{t_j}\big( E_h(s)&-E_h(t_j)\big)A_hP_hv\,\mathrm{d} s\Big\|^2
=
\Big\|\sum_{j=1}^n \sum_{i=1}^{\mathcal{N}_h}\int_{t_{j-1}}^{t_j}\big(e^{-s\lambda_{i,h}^2}-e^{-t_j\lambda_{i,h}^2}\big)\lambda_{i,h} (P_hv,e_{i,h})e_{i,h}\dd s\,\Big\|^2
\\
&=
\sum_{i=1}^{\mathcal{N}_h}
\bigg[
\Big|\sum_{j=1}^n \int_{t_{j-1}}^{t_j}\lambda_{i,h}^2e^{-s\lambda_{i,h}^2}\lambda^{-\frac\varrho2}_{i,h}\big(1-e^{-(t_j-s)\lambda_{i,h}^2}\big)\,\dd s\,\Big|^2  \lambda^{\varrho-2}_{i,h}(P_hv,e_{i,h})^2
\bigg]
\\
&\leq
Ck^{\frac\varrho2} \sum_{i=1}^{\mathcal{N}_h}
\Big| \int_{0}^{t_n}\lambda_{i,h}^2e^{-s\lambda_{i,h}^2}\,\dd s\,\Big|^2  \lambda^{\varrho-2}_{i,h}(P_hv,e_{i,h})^2
\\
&\leq
Ck^{\frac\varrho2} \sum_{i=1}^{\mathcal{N}_h}\lambda^{\varrho-2}_{i,h}(P_hv,e_{i,h})^2
\\
&=
Ck^{\frac\varrho2}\|A_h^{\frac{\varrho-2}2}P_hv\|^2
\\
& \leq
 Ck^{\frac\varrho2}|v|_{\varrho-2}^2.
\end{split}
\end{align}
Likewise,
\begin{align*}
\begin{split}
\Big\|k \sum_{j=1}^n ( E_h(t_j)-E_{k,h}^j)A_hP_hv\,\Big\|^2
&=
\Big\|\sum_{j=1}^n k \sum_{i=1}^{\mathcal{N}_h}(e^{-t_j\lambda_{i,h}^2}-r(k \lambda_{i,h}^2)^j)\lambda_{i,h} (P_hv,e_{i,h})e_{i,h}\,\Big\|^2
\\
&=
\sum_{i=1}^{\mathcal{N}_h}
\Big|\sum_{j=1}^nk (e^{-t_j\lambda_{i,h}^2}-r(k \lambda_{i,h}^2)^j)\lambda_{i,h}\,\Big|^2  (P_hv,e_{i,h})^2.
\end{split}
\end{align*}
Here we consider two possibilities: either $k\lambda^2_{i,h} \leq 1$ or $k\lambda^2_{i,h} > 1$.
For all summands with $k \lambda_{i,h}^2\leq 1$,  we rely on \eqref{eq:ez-r(z)-error} to get
\begin{align*}
\begin{split}
\Big|k\sum_{j=1}^n\Big(e^{-jk \lambda_{i,h}^2} -r(k \lambda_{i,h}^2)^j\Big) \lambda_{i,h}\Big|
&
\leq
C\lambda_{i,h}^5  k^2\sum_{j=1}^njk e^{-c(j-1)k \lambda_{i,h}^2}
\leq
C\lambda_{i,h}^5  k\int_0^{\infty} (r+k)   e^{-c r \lambda_{i,h}^2}\,\dd r
\\
&
\leq
C\lambda_{i,h}^5  k\Big( \frac1{(c \lambda_{i,h}^2)^2}+\frac k {c\lambda_{i,h}^2}\Big)
\leq
C \lambda_{i,h}^{\frac{\varrho-2}2} k^{\frac\varrho4}.
\end{split}
\end{align*}
For all summands with $k \lambda_{i,h}^2>1$, utilizing the fact $\sup_{s \in[0,\infty)} se^{-s} <\infty$ yields
\begin{align*}
\begin{split}
\Big|k\sum_{j=1}^n\Big(e^{-jk \lambda_{i,h}^2} -r( k \lambda_{i,h}^2)^j\Big) \lambda_{i,h}\Big|
&\leq
C \Big(k \lambda_{i,h}e^{-k \lambda_{i,h}^2} \sum_{j=1}^n e^{-(j-1)}+ \frac{k \lambda_{i,h}}{1+ k\lambda^2_{i,h}}\sum_{j=1}^n 2^{-j+1}\Big)
\\
&\leq
Ck^{\frac\varrho4}  \lambda^{\frac{\varrho-2}2}_{i,h} (k\lambda^2_{i,h})^{1-\frac \varrho4}\Big( e^{-k \lambda_{i,h}^2} + \frac1 {k\lambda^2_{i,h}}\Big)
\\
&\leq
Ck^{\frac\varrho4}\lambda^{\frac{\varrho-2}2}_{i,h}.
\end{split}
\end{align*}
This together with \eqref{eq:bound-Ph-H2} and \eqref{eq:Ah-A-bound}
proves
\begin{align}\label{eq:estimate-Phi(k,h)(s)-integrand-III}
\Big\|k\sum_{j=1}^n\big( E_h(t_j)-E_{k,h}^j\big)A_hP_hv\,\Big\|^2
\leq
Ck^{\frac\varrho2}\sum_{i=1}^{\mathcal{N}_h}\lambda^{\varrho-2}_{i,h}(P_hv,e_{i,h})^2
\leq
Ck^{\frac\varrho2}\|A_h^{\frac{\varrho-2}2}P_hv\|^2
\leq
Ck^{\frac\varrho2}|v|_{\varrho-2}^2.
\end{align}
Finally, plugging \eqref{eq:estimate-Phi(k,h)(s)-integrand-I}-\eqref{eq:summand-first-full} and \eqref{eq:estimate-Phi(k,h)(s)-integrand-III} into \eqref{eq:estimate-Phi(k,h)(s)-integrand-decompose} shows \eqref{lem:error-deterministc-potial-full-integrand-II} and thus finishes the proof.

\bibliography{bibfile}

\bibliographystyle{siam}

\end{document}